\pdfoutput=1
\documentclass[10pt]{article}
\usepackage{graphicx}
\usepackage{amssymb}
\usepackage{epstopdf}
\usepackage{amsmath}
\usepackage{amsfonts}

\def\phi{{\varphi}}

\def\P{{\mathcal P}}

\DeclareSymbolFont{AMSb}{U}{msb}{m}{n}
\DeclareMathSymbol{\N}{\mathbin}{AMSb}{"4E}
\DeclareMathSymbol{\Z}{\mathbin}{AMSb}{"5A}
\DeclareMathSymbol{\R}{\mathbin}{AMSb}{"52}
\DeclareMathSymbol{\Q}{\mathbin}{AMSb}{"51}
\DeclareMathSymbol{\I}{\mathbin}{AMSb}{"49}
\DeclareMathSymbol{\C}{\mathbin}{AMSb}{"43}

\def\be{\begin{equation}}
\def\ber{\begin{eqnarray}}
\def\eer{\end{eqnarray}}

\def\vv{{\bf v}}

\def\beq{\begin{equation}}
\def\eeq{\end{equation}}
\newcommand{\E}{\mathbb{E}}
 \newcommand{\tens}{%
  \mathbin{\mathop{\otimes}}%
}

\begin{document}

\addtolength{\textheight}{0 cm} \addtolength{\hoffset}{0 cm}
\addtolength{\textwidth}{0 cm} \addtolength{\voffset}{0 cm}
\def\XX{\mathcal X}

\newcommand{\aaa}{\mathbb{A}}
\newcommand{\bb}{\mathbb{B}}
\newcommand{\cc}{\mathbb{C}}
\newcommand{\dd}{\mathbb{D}}
\newcommand{\ee}{\mathbb{E}}
\newcommand{\mm}{\mathbb{M}}
\newcommand{\rr}{\mathbb{R}}
\newcommand{\pp}{\mathbb{P}}
\newcommand{\qq}{\mathbb{Q}}
\newcommand{\ttt}{\mathbb{T}}
\newcommand{\zz}{\mathbb{Z}}

\def\vep{\varepsilon}
\def\<{\langle}
\def\>{\rangle}
\def\dsubset{\subset\subset}
\newcommand{\as}[1]{\begin{align*}#1\end{align*}}
\newcommand{\ald}[1]{\begin{aligned}#1\end{aligned}}

\newcommand{\Var}{\mathrm{Var}}
\newcommand{\mcal}[1]{\mathcal{#1}}

\def\AA{\mathcal A}
\def\BB{\mathcal B}
\def\CC{\mathcal C}
\def\DD{\mathcal D}
\def\FF{\mathcal F}
\def\GG{\mathcal G}
\def\EE{\mathcal E}
\def\JJ{\mathcal J}
\def\KK{\mathcal K}
\def\LL{\mathcal L}
\def\MM{\mathcal M}
\def\PP{\mathcal P}
\def\SS{\mathcal S}
\def\VV{\mathcal V}
\def\T{\mathcal T}
\def\d{\, {\rm d}}
 \newcommand{\lf}{\left}
\newcommand{\rt}{\right}
\newcommand{\vphi}{\varphi}
\newcommand{\no}{\nonumber}
\newcommand{\eq}[1]{\begin{equation}#1\end{equation}}
\newcommand{\eqs}[1]{\begin{equation*}#1\end{equation*}}

\newcommand{\ZZ}{\mathbb{Z}}
\newcommand{\Rm}{\mathbb{R}}
\newcommand{\RR}{\mathbb{R}}
\newcommand{\NN}{\mathbb{N}}
\newcommand{\sU}{\mathcal{U}}
\newcommand{\sF}{\mathcal{F}}
\newcommand{\sM}{\mathcal{M}}
\newcommand{\sS}{\mathcal{S}}
\newcommand{\mL}{\mathcal{L}}
\newcommand{\ac}{\hbox{\small ac}}
\newcommand{\mC}{\ensuremath{\mathcal{C}}}
\newcommand{\mU}{\ensuremath{\mathcal{U}}}
\newcommand{\mT}{\ensuremath{\mathcal{T}}}
\newcommand{\mS}{\ensuremath{\mathcal{S}}}
\newcommand{\mF}{\ensuremath{\mathcal{F}}}
\newcommand{\Nm}{\ensuremath{\mathbb{N}}}
\newcommand{\Zm}{\ensuremath{\mathbb{Z}}}
\newcommand{\Hm}{\ensuremath{\mathbb{H}}}
\newcommand{\mM}{\ensuremath{\mathcal{M}}}
\newcommand{\mK}{\ensuremath{\mathcal{K}}}
\newcommand{\mD}{\ensuremath{\mathcal{D}}}
\newcommand{\mA}{\ensuremath{\mathcal{A}}}
\newcommand{\mO}{\ensuremath{\mathcal{O}}}
\newcommand{\mI}{\ensuremath{\mathcal{I}}}
\newcommand{\mB}{\ensuremath{\mathcal{B}}}
\newcommand{\Tm}{\ensuremath{\mathbb{T}}}
\newcommand{\mE}{\ensuremath{\mathcal{E}}}
\newcommand{\vs}{\vspace{.5cm}}

\def\proof {\noindent{\sc{Proof. }}}
\def\qed {\mbox{}\hfill {\small \fbox{}} \\}  
\def\lto{\longrightarrow}
\def\lmto{\longmapsto}
\def\eq{\Longleftrightarrow}
\def\leq{\leqslant}
\def\geq{\geqslant}
\def \uK {K^+}
\def \oK {K^-}
\def \calT {\mathcal T}

\newtheorem{lem}{Lemma}
\newtheorem{thm}{Theorem}
\newtheorem{conj}[lem]{Conjecture}
\newtheorem{ques}[lem]{Question}
\newtheorem{cor}[lem]{Corollary}
\newtheorem{prop}[lem]{Proposition}
\newtheorem{defn}[lem]{Definition}
\newtheorem{note}[lem]{Note}
\newtheorem{rmk}{Remark}
\def\proof {\noindent{\sc{Proof. }}}
\def\qed {\mbox{}\hfill {\small \fbox{}} \\}  

\newcommand{\grad}{\operatorname{grad}}
\newcommand{\Leg}{\mathcal{L}}

\newcommand{\lbstoc}[0]{\underline{B}}
\newcommand{\ubstoc}[0]{\overline{B}^{\text{stoc}}}
\newcommand{\deriv}[2]{\ensuremath{\frac{d{#1}}{d{#2}}}}
\newcommand{\Id}[1]{\ensuremath{\boldsymbol{1}[#1]}}
\newcommand{\abs}[1]{\ensuremath{\left\lvert#1\right\rvert}}
\renewcommand{\P}[1]{\ensuremath{\mathbb{P}(#1)}}
\newcommand{\bracket}[1]{\ensuremath{\left(#1\right)}}
\newcommand{\expect}[2][]{\ensuremath{\mathbb{E}_{#1}\left[#2\right]}}
\newcommand{\expcond}[2]{\ensuremath{\mathbb{E}\left[#1\middle\rvert #2\right]}}
\newcommand{\ball}[2]{\ensuremath{B(#1,#2)}}
\newcommand{\proj}[1]{\ensuremath{\text{proj}_{#1}}}
\newcommand{\pderiv}[2]{\ensuremath{\frac{\partial{#1}}{\partial{#2}}}}
\newcommand{\Cmik}[0]{\ensuremath{C}}

\title{A Theory of Transfers: Duality and convolution}
\author{Malcolm Bowles\thanks{This is part of the PhD dissertation of this author at the University of British Columbia.} \quad   and \quad Nassif  Ghoussoub\thanks{Partially supported by a grant from the Natural Sciences and Engineering Research Council of Canada. } 
\\ \\
{\it\small Department of Mathematics,  University of British Columbia}\\
{\it\small Vancouver BC Canada V6T 1Z2}\\
}

\date{April 16, 2018 (revised on October 24, 2018)}
\maketitle

\begin{abstract} We introduce and study the permanence properties of {\it the class of linear transfers} between probability measures and {\it the dual class of Kantorovich operators} between continuous functions. The class of linear transfers contains all cost minimizing mass transports, but also {\it Balayage operations}, {\it martingale mass transports}, the {\it Schr\"odinger bridge} associated to a reversible Markov process, optimal Skorokhod embeddings, and  the {\it weak mass transports} of Talagrand, Marton, Gozlan and others. The class also includes various stochastic mass transports to which Monge-Kantorovich theory does not apply. We also introduce {\it the cone of convex transfers}, which include any $p$-power ($ p \geq 1$) of a linear transfer, but also the logarithmic entropy, optimal mean field plans, the Donsker-Varadhan information, and certain free energy functionals. This first paper is mostly focused on  exhibiting examples that point to the pervasiveness of the concept in the important quest of  correlating probability distributions. Duality formulae for general transfer inequalities follow in a very natural way. We also associate to each linear transfer, a corresponding {\it effective transfer} (or a {\em generalized Peierls barrier}) and a dual {\it effective Kantorovich operator}, that could be seen as a generalization of the effective Lagrangians and Hamiltonians in weak KAM theory. In a forthcoming paper, we show how it  allows, in particular, for the development of a stochastic counterpart of the Fathi-Mather theory.

\end{abstract} 

\tableofcontents

\section{Introduction}
Stochastic control problems and several other analytical and statistical procedures that correlate two probability distributions share many of the useful properties of optimal mass transportation between probability measures. However, these correlations often lack at least two of the useful features of Monge-Kantorovich theory \cite{V}. For one, they are not symmetric, meaning that the problem imposes a specific direction from one of the marginal distributions to the other. Moreover, many of those do not arise as cost minimizing problems associated to functionals $c(x, y)$ that assign ``a price for moving one particle $x$ to another $y$." As such, they are not readily amenable to the duality theory of Monge-Kantorovich. In this paper, we isolate and study a notion of transfers between probability measures that encapsulates both the deterministic and stochastic versions of transport problems studied by Mikami-Thieulin \cite{M-T}) and Barton-Ghoussoub \cite{B-G}, the optimal Skorokhod embeddings of Ghoussoub-Kim-Pallmer \cite{G-K-P3}, but also includes the {\it weak mass transports} of Talagrand \cite{Ta1, Ta2}, Marton \cite{Ma1, Ma2} and Gozlan et al. \cite{GL, Go4}, the logarithmic entropy, optimal mean field plans \cite{Sav}, the Donsker-Varadhan information \cite{DV}, and many other energy correlation functionals.

This first paper introduces the unifying concepts of {\it linear and convex mass transfers} and exhibits several examples that illustrate the potential scope of this approach. The underlying idea has been implicit in many related works and should be familiar to the experts. But, as we shall see, the systematic study of these structures add clarity and understanding, allow for non-trivial extensions, and open up a whole new set of interesting problems. The ultimate purpose is to extend many of the remarkable properties enjoyed by standard mass transportations to linear and convex transfers, and hence to the stochastic case, or at least to weak mass transports. This is a vast undertaking. We therefore decided --for this first paper-- to give a sample of the results that can be inspired and eventually extended from the standard theory of mass transport.  We  chose to state here the most basic permanence properties of the cones of transfers and to establish general duality formulas for potential comparisons between different transfers that extend the work of Bobkov-G\"{o}tze \cite{B-G}, Gozlan-Leonard \cite{GL}, Maurey \cite{Mau} and others. We also show how the approach of Bernard-Buffoni \cite{B-B1, B-B2} on Fathi's  weak KAM \cite{Fa} and Mather theory \cite{Mat} extend to linear transfers, and therefore could, for example,  be applied to the stochastic case. We will pursue this in a forthcoming paper \cite{B-G1}. Furthermore, we shall present in \cite{B-G2} a notion of {\it linear and convex multi-transfers} between several probability distributions that will --among other things-- extend the theory of multi-marginal mass transportation.   
 
We shall focus here on probability measures on compact spaces, even though the right settings for most applications and examples are complete metric spaces, or at least $\R^n$. This will allow us to avoid the usual functional analytic complications, and concentrate on the algebraic aspects of the theory. The simple compact case will at least point to results that can be expected to hold and be proved --albeit with additional analysis and suitable hypothesis -- in more general situations. In the case of $\R^n$, which is the setting for many examples stated below, the right duality is between the space $Lip (\R^n)$ of all bounded and Lipschitz functions and the space of Radon measures with finite first moment. 

 With this in mind, we shall denote by $C(X)$ (resp. $USC(X)$), (resp $LSC(X)$) to be the spaces of continuous (resp., upper semi-continuous), (resp., lower semi-continuous) functions on a compact space $X$. The class of signed (resp., probability) measures  on $X$ will be denoted by ${\cal M}(X)$ (resp., ${\cal P}(X)$). 

If now ${\mathcal T}: {\mathcal M}(X)\times {\mathcal M}(Y)  \to \R\cup \{+\infty\}$ is a proper convex functional, we shall denote by $D(\calT)$ its effective domain, that is the set where it takes finite values. We shall always assume that $D(\calT)\subset {\mathcal P}(X)\times {\mathcal P}(Y)$, where ${\mathcal P}(X)$ is the set of probability measures on $X$. The ``partial domains" of ${\mathcal T}$ are then denoted by,
$$D_1(\calT)=\{\mu \in {\mathcal P}(X); \exists \nu\in {\mathcal P}(Y), (\mu, \nu)\in D(\calT)\}
\hbox{
\,\, and  \,\,
$D_2(\calT)=\{\nu \in {\mathcal P}(Y); \exists \mu\in {\mathcal P}(X), (\mu, \nu)\in D(\calT)\}.$}
$$
 We consider for each $\mu \in {\mathcal P}(X)$ (resp.,  $\nu \in {\mathcal P}(Y)$) the partial maps ${\mathcal T}_\mu$ on ${\mathcal P}(Y)$ (resp., ${\mathcal T}_\nu$ on ${\mathcal P}(X)$) given by $\nu \to {\mathcal T} (\mu, \nu)$ (resp., $\mu \to {\mathcal T} (\mu, \nu)$). 
 
  \begin{defn} \label{notion1} \rm Let $X$ and $Y$ be two compact spaces, and let 
${\mathcal T}: {\mathcal P}(X)\times {\mathcal P}(Y) \to \R\cup \{+\infty\}$ be a proper bounded below, convex and weak$^*$ lower semi-continuous functional on ${\mathcal M}(X)\times {\mathcal M}(Y)$. We say that 
\begin{enumerate}
\item $\calT$ is a {\it backward linear transfer}, if  there exists a map $T ^-: C(Y) \to LSC(X)$ such that for each $\mu \in D_1(\calT)$, the Legendre transform of ${\mathcal T}_\mu$ on ${\mathcal M}(Y)$  satisfies:
\begin{equation}
\hbox{${\mathcal T}^*_\mu (g)=\int_XT ^-g(x) \, d\mu(x)$ \quad for any $g\in C(Y)$}.
\end{equation}

\item $\calT$ is a {\it forward linear transfer}, if there exists a map $T ^+: C(X)\to USC(Y)$ such that for each $\nu \in D_2(\calT)$, the Legendre transform of  ${\mathcal T}_\nu$ on ${\mathcal M}(X)$ satisfies:
 \begin{equation}
\hbox{${\mathcal T}^*_\nu (f)=-\int_Y {T ^+}(-f)(y) \, d\nu(y)$ \quad for any $f\in C(X)$. }
\end{equation}  
  \end{enumerate}
 We shall call $T ^+$ (resp., $T ^-$) the {\it forward (resp., backward) Kantorovich operator} associated to ${\mathcal T}$.
   \end{defn}
By Legendre transform of ${\mathcal T}_\nu$, we mean here  
\[
{\mathcal T}^*_\nu (f)=\sup\{\int_Xf d\mu -{\mathcal T}_\nu(\mu);\mu \in {\mathcal P}(X) \}=\sup\{\int_X f d\mu -{\mathcal T}(\mu, \nu); \, \mu\in {\mathcal P}(X) \}.
\]
This is because we are assuming that ${\mathcal T}_\nu$ and ${\mathcal T}_\mu$ are equal to $+\infty$ whenever $\mu$ and $\nu$ are not probability measures.  
So, if ${\mathcal T}$ is a forward linear transfer on $X\times Y$, then for any  $\mu \in {\mathcal P}(X)$ and $\nu \in {\mathcal P}(Y)$, we have
 \begin{equation}
{\mathcal T}(\mu, \nu)= \sup\big\{\int_{Y}{T ^+}f(y)\, d\nu(y)-\int_{X}f(x)\, d\mu(x);\,  f \in C(X) \big\},
\end{equation}
while if ${\mathcal T}$ is a backward linear transfer on $X\times Y$, then 
 \begin{equation}
{\mathcal T}(\mu, \nu)=\sup\big\{\int_{Y}g(y)\, d\nu(y)-\int_{X}{T ^-}g(x)\, d\mu(x);\,  g \in C(Y)\big\}.
\end{equation}
We shall say that {\it a transfer ${\cal T}$ is symmetric} if 
$$\hbox{${\cal T}(\nu, \mu):={\cal T}(\nu, \mu)$ for all $\mu \in {\mathcal P}(X)$ and $\nu \in {\mathcal P}(Y)$}.
$$
Note that if ${\cal T}$ is a backward linear transfer with Kantorovich operator $T^-$, then $\tilde {\cal T}(\mu, \nu):={\cal T}(\nu, \mu)$ is a forward linear transfer with Kantorovich operator ${\tilde T}^+f=-T^-(-f)$. This means that if  ${\cal T}$ is symmetric, then $T^+f=-T^-(-f)$.\\
The class of linear transfers is quite large and ubiquitous in analysis. To start with, it contains all cost minimizing mass transports, that is functionals on ${\mathcal P}(X)\times {\mathcal P}(Y)$ of the form,
\begin{eqnarray}
{\mathcal T}_c(\mu, \nu):=\inf\big\{\int_{X\times Y} c(x, y)) \, d\pi; \pi\in \mK(\mu,\nu)\big\},
\end{eqnarray}
where $c(x, y)$ is a continuous cost function on the product measure space $X\times Y$, and $\mK(\mu,\nu)$ is the set of probability measures $\pi$ on $X\times Y$ whose marginal on $X$ (resp. on $Y$) is $\mu$ (resp., $\nu$) {\it (i.e., the transport plans)}. A consequence of the Monge-Kantorovich theory is that cost minimizing transports ${\mathcal T}_c$ are both forward and  backward linear transfers. The {\it Schr\"odinger bridge} problem associated to a reversible Markov process \cite{GLR} is also a symmetric backward and forward linear transfer.\\
Other examples, which are only one-directional linear transfers, are  the various {\it Martingale mass transports}, the {\it weak mass transports} of Marton, Gozlan and collaborators. However, what motivated us to develop the concept of transfers are the stochastic mass transports, which do not minimize a given cost function between point particles, since the cost of transporting a Dirac measure to another is often infinite. This said, we should show however that if the set $\{\delta_x; x\in X\}$ is contained in $D_1({\mathcal T})$, then we can represent such a linear transfer as a generalized mass transport, a notion recently formalized by Gozlan et al. \cite{Go4}. 

Note that we did not specify any property on the maps $T^+$ and $T^-$. However, the fact that they arise from a Legendre transform imposes on them certain properties such as those exhibited in the following. 

\begin{defn}
If $X$ and $Y$ are two compact spaces, say that a map $T ^-: C(Y) \to LSC(X)$ (resp., $T ^+: C(X)\to USC(Y)$) is  {\it a convex operator} (resp., {\it a concave operator}), if it satisfies the following conditions:
 \begin{enumerate}
 \item If $f_1\leq f_2$ in $C(Y)$ (resp., in $C(X)$), then $T^-f_1\leq T^-f_2$ (resp., $T^+f_1\leq T^+f_2$).
 \item For any $\lambda \in [0, 1]$, $f_1, f_2$ in $C(Y)$ (resp., in $C(X)$), we have 
\begin{equation*}
\hbox{$T^-(\lambda f_1+(1-\lambda)f_2)\leq \lambda T^-f_1+(1-\lambda)T^-f_2$ (resp., $T^+(\lambda f_1+(1-\lambda)f_2)\geq \lambda T^+f_1+(1-\lambda)T^+f_2$.}
\end{equation*}
\item For any constant $c\in \R$ and $f\in C(Y)$ (resp., $C(X)$), there holds that $T^-(f+c)=T^-f +c$ (resp., $T^+(f+c)=T^+f +c$.
\item $T^-$ (resp., $T^+$) is $1-$Lipshitz, i.e. $\|T^-f_1- T^-f_2\|\leq \|f_1-f_2\|$. 
\item If $(f_n)_n, f$ in $C(Y)$ (resp., in $C(X)$) are such that $f_n \to f$ weakly, then 
\begin{equation*}
 T^-f \leq \liminf_nT^-f_n \quad \hbox{(resp., $T^+f \geq \limsup_nT^-f_n$).}
\end{equation*}
\end{enumerate}
Note that $T^-$ (resp., $T^+$) extend --with the same properties-- to operators $T^-:LSC (Y)\to LSC (X)$ (resp., $T^+: USC(X)\to LSC (Y)$). 
\end{defn}
We leave it to the reader to check that the backward (resp., forward) maps in Definition \ref{notion1} are necessarily convex (resp., concave) operators. 

Note also that conversely, any convex operator $T$ (resp., concave) defines a backward (resp., forward) linear transfer via the formula
 \begin{equation}
{\mathcal T}(\mu, \nu)=\left\{ \begin{array}{llll}
\sup\big\{\int_{Y}g(y)\, d\nu(y)-\int_{X}{T}g(x)\, d\mu(x);\,  g \in C(Y)\big\} \quad &\hbox{if $\mu \in {\mathcal P}(X),  \nu \in {\mathcal P}(Y)$,}\\
+\infty \quad &\hbox{\rm otherwise.}
\end{array} \right.
\end{equation}
Indeed, it is clear that ${\mathcal T}_\mu  \geq \Gamma_{_{T, \mu}}^*$, where $\Gamma_{_{T, \mu}}$ is the convex continuous function on $C(Y)$ defined by $\Gamma_{_{T, \mu}}(g)=\int_{X}{T }g(x)\, d\mu(x)$ and that ${\mathcal T}_\mu  = \Gamma_{_{T, \mu}}^*$ on the probability measures on $Y$.  If now $\nu$ is a positive measure with $\lambda:=\nu (Y) >1$, then 
\[
\Gamma_{_{T, \mu}}^*(\nu)=\sup\big\{\int_{Y}g(y)\, d\nu(y)-\int_{X}{T }g(x)\, d\mu(x);\,  g \in C(Y)\big\}\geq n\lambda -\int_XT(n)\d\mu=n(\lambda -1)-\int_XT(0)\d\mu,
\]
where we have used property (3) to say that $T(n)=n+T(0)$. Hence $\Gamma_{_{T, \mu}}^*(\nu)=+\infty$. A similar reasoning applies to when $\lambda<1$ and for the concave case.

The class of  linear transfers has remarkable permanence properties. The two most important ones are stability under inf-convolution and tensorization, which allow to create an even richer class of transfers, such as {\it the ballistic stochastic optimal transport} and  {\it broken geodesics of transfers.}  
However, a natural and an even richer family of transfers is the class of {\it convex transfers,} which are essentially suprema of linear transfers.  
\begin{defn} \rm 
A proper convex and weak$^*$ lower semi-continuous functional ${\mathcal T}: {\mathcal P}(X)\times {\mathcal P}(Y) \to \R\cup \{+\infty\}$ is said to be a {\it backward convex transfer} (resp., {\it forward convex transfer}), if there exists a family of backward linear transfers (resp., forward linear transfers) $({\mathcal T}_i)_{i\in I}$ such that for all $\mu\in {\cal P}(X)$, $\nu \in {\cal P}(Y)$, 
\begin{equation}
{\mathcal T}(\mu, \nu)=\sup_{i\in I}{\mathcal T}_i(\mu, \nu).
\end{equation}
In other words, a backward convex transfer (resp., forward convex transfer) can be written as:
 \begin{equation}
{\mathcal T}(\mu, \nu)=\sup\big\{\int_{Y}g(y)\, d\nu(y)-\int_{X}{T_i^-}g(x)\, d\mu(x);\,  g \in C(Y), i\in I\big\},
\end{equation}
respectively,
  \begin{equation}
{\mathcal T}(\mu, \nu)= 
\sup\big\{\int_{Y}{T_i^+}f(y)\, d\nu(y)-\int_{X}f(x)\, d\mu(x);\,  f \in C(X), i\in I \big\}.
\end{equation}
where $(T_i ^-)_{i\in I}$ (resp., $(T_i ^+)_{i\in I}$) is a family of convex operators from $C(Y) \to LSC(X)$ (resp., concave operators from $C(X) \to LSC(Y)$).  
 \end{defn}
 In addition to linear transfers, we shall see that any $p$-power ($ p \geq 1$) of a linear transfer is a convex transfer in the same direction.   More generally, for any convex increasing real function $\gamma$ on $\R^+$ and any linear backward (resp., forward) transfer, the map $\gamma ({\cal T})$ is a convex backward (resp., forward) transfer. 

Note that if a ${\cal T}$ is convex backward (resp., forward) transfer, then   
\begin{equation}
\hbox{${\mathcal T}_\mu=(S_\mu^-)^*$\quad  and \quad ${\mathcal T}_\nu=(S_\nu^+)^*$},
\end{equation}
where $S_\mu^-(g)=\inf\limits_{i\in I}\int_XT_i ^-g(x) \, d\mu(x)$ for $g\in C(Y)$ and $S_\nu^+(f)=\sup\limits_{i\in I}\int_Y {T_i ^+}(-f)(y) \, d\nu(y)$ \, for $f\in C(X)$. However, we only have 
\begin{equation}
\hbox{${\mathcal T}^*_\mu\leq S_\mu^-$\quad  and \quad ${\mathcal T}^*_\nu\leq -S_\nu^+$}, 
\end{equation}
since $S_\mu^-$ (resp., $S_\nu^+$) are not necessarily convex (resp., concave). 
We can therefore introduce the notions of {\it completely convex} transfers for when we have equality above, that is when  $S_\mu^-$ is a convex operator (resp., $S_\mu^+$ is concave) and ${\mathcal T}^*_\mu = S_\mu^-$ (resp., ${\mathcal T}^*_\nu= -S_\nu^+$).  For instance, this will be the case for the following {\it generalized entropy,}
 \begin{equation}
 {\cal T}(\mu, \nu)=\int_X \alpha (\frac{d\nu}{d\mu})\,  d\mu, \quad \hbox{if $\nu<<\mu$ and $+\infty$ otherwise,}
 \end{equation}
which is a backward completely convex transfer, whenever $\alpha$ is a strictly convex lower semi-continuous superlinear real-valued function on $\R^+$. 
 The important example of the {\it logarithmic entropy} 
 \begin{equation}
 {\cal H}(\mu, \nu)=\int_X \log (\frac{d\nu}{d\mu})\,  d\nu, \quad \hbox{if $\nu<<\mu$ and $+\infty$ otherwise,}
\end{equation}
is of course one of them, but it is much more as we  now focus on  a remarkable subset of the cone of completely convex transfers, which is the class of {\it entropic transfers}, that we define as follows:

\begin{defn} \rm  Let $\alpha$ (resp., $\beta$) be a convex increasing (resp., concave increasing) real function on $\R$, and let ${\mathcal T}: {\mathcal P}(X)\times {\mathcal P}(Y) \to \R\cup \{+\infty\}$ be a proper (jointly) convex and weak$^*$ lower semi-continuous functional. We say that 
\begin{itemize}
\item ${\mathcal T}$ is a {\it $\beta$-backward transfer}, if there exists a convex operator $T ^-: C(Y) \to LSC(X)$ such that for each $\mu \in D_1(\calT)$, the Legendre transform of ${\mathcal T}_\mu$ on ${\mathcal M}(Y)$ satisfies:
$$\hbox{${\mathcal T}^*_\mu (g)=\beta \left(\int_XT ^-g(x) \, d\mu(x)\right)$ \quad for any $g\in C(Y)$}.
$$

\item $\calT$ is a {\it $\alpha$-forward transfer}, if there exists a concave operator $T ^+: C(X)\to USC(Y)$ such that for each $\nu \in D_2(\calT)$, the Legendre transform of  ${\mathcal T}_\nu$ on ${\mathcal M}(X)$ satisfies: 
$$\hbox{${\mathcal T}^*_\nu (f)=-\alpha \left(\int_Y {T ^+}(-f)(y) \, d\nu(y)\right)$ \quad for any $f\in C(X)$. }
$$  
  \end{itemize}
 \end{defn}
So, if ${\mathcal T}$ is an $\alpha$-forward transfer on $X\times Y$, then for any probability measures $\mu \in {\mathcal P}(X)$ and $\nu \in {\mathcal P}(Y)$, we have
 \begin{equation}
{\mathcal T}(\mu, \nu)= 
\sup\big\{\alpha\left(\int_{Y}{T ^+}f(y)\, d\nu(y)\right)-\int_{X}f(x)\, d\mu(x);\,  f \in C(X) \big\},
\end{equation}
while if ${\mathcal T}$ is a $\beta$-backward transfer, then 
 \begin{equation}
{\mathcal T}(\mu, \nu)=\sup\big\{\int_{Y}g(y)\, d\nu(y)-\beta \left(\int_{X}{T ^-}g(x)\, d\mu(x)\right);\,  g \in C(Y)\big\}.
\end{equation}
Entropic transfers are completely convex transfers. A typical example is  of course the logarithmic entropy, since it can be written as 
\begin{equation}
 {\cal H}(\mu, \nu)=\sup\{\int_X f\, d\nu-\log(\int_Xe^{f}\, d\mu);\,  f\in C(X)\}, 
\end{equation}
making it a $\log$-backward transfer. The {\em Donsker-Varadhan information} is defined 
as 
\begin{equation}
{\cal I}(\mu, \nu):=\begin{cases}\EE(\sqrt{f}, \sqrt{f}), \ \ &\text{ if }\ \mu=f\nu, \sqrt{f}\in\dd(\EE)\\
+\infty, &\text{ otherwise,}
\end{cases}
\end{equation}
where $\EE$ is a Dirichlet form  with domain $\dd(\EE)$ on
$L^2(\nu)$. It  is another example of a backward completely convex transfer, since it can also be written as
\begin{equation}
{\cal I}(\mu, \nu)=\sup\{\int_X f\, d\nu-\log \|P_1^f\|_{L^2(\mu)};\,  f\in C(X)\}, 
\end{equation}
 where $P_t^f$ is an associated (Feynman-Kac) semi-group of operators on $L^2(\nu)$. 
 More examples of $\alpha$-forward transfers and $\beta$-backward transfers with readily computable Kantorovich operators can be obtained by convolving entropic transfers  with linear transfers of the same direction. 
 
 In section 7, we show how the concepts of linear and convex transfers lead naturally to more transparent proofs and vast extensions, of many well known duality formulae for transport-entropy inequalities. 
 
   In section 8, we associate to any backward linear transfer ${\mathcal T}$ and its Kantorovich operator $T$, a corresponding {\it effective backward linear transfer} ${\mathcal T}_\infty$ and its {\it effective Kantorovich operator} $T_\infty$. ${\mathcal T}_\infty$ is obtained by an infinite  inf-convolution process, while $T_\infty$ is obtained by an infinite iteration procedure, which lead to fixed points for such a non-linear operator. In the case of a linear Kantorovich operator (the push-forward transfer) this reduces to classical ergodic theory.  If now ${\mathcal T}$ is the optimal mass transport minimizing a cost given by the generating function of a Lagrangian $L$ on a compact manifold $M$, that is 
 \begin{equation}\label{GF}
c^L(y,x):=\inf\{\int_0^1L(t, \gamma(t), {\dot \gamma}(t))\, dt; \gamma\in C^1([0, 1), M);  \gamma(0)=y, \gamma(1)=x\},
\end{equation}
 the Kantorovich operator is then given by the Lax-Oleinik semi-group, whose fixed points correspond to weak KAM solutions as described by Fathi \cite{Fa}. The extension of this result to general transfers allows for a similar approach for the stochastic counterpart of Mather theory. This will be the subject of a companion paper \cite{B-G1}. 
 
 \section{Examples of linear mass transfers}

\subsection{Elementary linear mass transfers}

    \noindent {\it Example 1: The push-forward transfer} 
    
 {\it The identity transfer} is the map ${\cal I}$ on ${\cal P}(X)\times {\cal P}(X)$ defined by 
 \begin{equation}
{\mathcal I}(\mu, \nu)=\left\{ \begin{array}{llll}
0 \quad &\hbox{if $\mu=\nu $}\\
+\infty \quad &\hbox{\rm otherwise.}
\end{array} \right.
\end{equation}
The corresponding Kantorovich operators are the identity map, that is $T^+f=T^-f=f$.
 
  More generally, if $\sigma$ is a continuous map from $X$ to $Y$, then 
   \begin{equation}
{\mathcal I}_\sigma (\mu, \nu)=\left\{ \begin{array}{llll}
0 \quad &\hbox{if $\sigma_\#\mu=\nu $}\\
+\infty \quad &\hbox{\rm otherwise.}
\end{array} \right.
\end{equation}
is a backward linear transfer with Kantorovich operator given by $T^-f=f\circ \sigma$. 

Similarly, any probability measure $\pi$ on $X\times Y$ induces a forward and backward linear transfer in the following way:
\begin{equation} \label{plan}
{\mathcal I}_\pi (\mu, \nu)=\left\{ \begin{array}{llll}
0 \quad &\hbox{if $\mu=\pi_1$ and $\nu=\pi_2.$}\\
+\infty \quad &\hbox{\rm otherwise,}
\end{array} \right.
\end{equation}
where $\pi_1$ (resp., $\pi_2$) is the first (resp., second) marginal of $\pi$. In this case,
\begin{equation}
T^-f(x)=\int_Yf(y)d\, \pi_x(y) \quad \hbox{and \quad $T^+f(y)=\int_Xf(x)d\, \pi_y(x)$,}
\end{equation}
where $(\pi_x)_x$ (resp., $(\pi_y)_y$) is the disintegration of $\pi$ with respect to $\pi_1$ (resp.,  $\pi_2$).\\

 \noindent {\it Example 2: The prescribed Balayage transfer}
 
 Given a convex cone of continuous functions ${\mathcal A} \subset C(X)$, where $X$ is a compact space, one can define an order relation between probability measures $\mu, \nu$ on $X$, called the  ${\mathcal A}$-balayage, in the following way. 
 \[
 \mu \prec_{\mathcal A} \nu \quad \hbox{ if and only if \quad $\int_X\phi \, d\mu \leq \int_X\phi \, d\nu$ for all $\phi$ in  ${\mathcal A}$.}
 \]
 Say that a probability measure $\pi$ on $X\times X$ is {\it an ${\cal A}$-dilation} if $\delta_x \prec_{\cal A} \pi_x$, where $(\pi_x)_x$ is the disintegration of $\pi$ with respect to its first marginal $\pi_1$. \\
To each ${\cal A}$-dilation $\pi$, one can define a backward linear transfer as in (\ref{plan}) above.
In this case, the corresponding backward Kantorovich transfer is again,
$T^-f(x)=\int_Yf(y)d\, \pi_x(y)$. \\

\noindent {\it Example 3: The prescribed Skorokhod transfer}
 
  Writing $Z\sim \rho$ if $Z$ is a random variable with distribution $\rho$, and letting $(B_t)_t$ denote Brownian motion, and ${\mathcal S}$  the corresponding class of --possibly randomized-- stopping times.  
   For a fixed $\tau\in {\mathcal S}$, one can associate a backward linear transfer in the following way:
  \begin{equation}
{\mathcal T}_{\tau}(\mu, \nu)=\left\{ \begin{array}{llll}
0 \quad &\hbox{if $B_0\sim \mu$ and $B_\tau \sim\nu$.}\\
+\infty \quad &\hbox{\rm otherwise.}
\end{array} \right.
\end{equation}  
Its backward Kantorovich operator is then $T^-f(x)=\mathbb{E}^{x}[f(B_\tau)]$, where the expectation is with respect to Brownian motion satisfying $B_0=x$.

\subsection{Zero-cost linear mass transfers}

 \noindent {\it Example 4: The null transfer} 
  
  This is simply the map ${\mathcal N}(\mu, \nu)=0$ for all probability measures $\mu$ on $X$ and $\nu$ on $Y$. 
  It is easy to see that it is both a backward and forward linear transfer with Kantorovich operators, 
 \begin{equation}
 \hbox{$T^-f\equiv \sup_{y\in Y}f(y)$ \quad and \quad $T^+f\equiv\inf_{x\in X}f(x)$.}
 \end{equation}
 Note that 
 \begin{equation}
 {\mathcal N}(\mu, \nu)=\inf\{{\mathcal I}_\sigma (\mu, \nu); \sigma: X\to Y\}=\inf\{{\mathcal I}_\pi (\mu, \nu); \pi \hbox{ is a transfer plan on $X\times Y$}\},
 \end{equation}
 where ${\mathcal I}_\sigma$ and ${\mathcal I}_\pi$ are the push-forward transfers defined in Example 1. This is a particular case, i.e., when the cost is trivial,  of a relaxation result of Kantorovich (e.g., see Villani \cite{V}).\\
 
  \noindent {\it Example 5: The Balayage transfer} 
  
 Let ${\mathcal A}$ be a closed convex cone in $C(X)$, and define now the {\em backward balayage transfer} ${\cal B}_b$ on ${\cal P}(X)\times {\cal P}(X)$ via
 \begin{equation}
{\mathcal B}^b(\mu, \nu)=\left\{ \begin{array}{llll}
0 \quad &\hbox{if $\mu \prec_{\mathcal A} \nu$}\\
+\infty \quad &\hbox{\rm otherwise.}
\end{array} \right.
\end{equation} 
A generalized version of a Theorem of Strassen \cite{St} yields the following relationship:

\begin{prop} Assume the cone ${\mathcal A}$ is proper, separates the points of $X$, and that it is stable under finite suprema. Then, for any two probability measures $\mu, \nu$ on $X$, the following are equivalent:
\begin{enumerate}
\item $\mu \prec_{\mathcal A} \nu$.
\item There exists an $\cal A$-dilation $\pi$ on $X\times X$ such that $\mu=\pi_1$ and $\nu=\pi_2$.
\end{enumerate}
\end{prop}
From this follows that 
\begin{equation}
{\mathcal B}^b(\mu, \nu)=\inf\{{\mathcal B}_\pi(\mu, \nu);  \hbox{$\pi$ is an $\cal A$-dilation}\}.
\end{equation}
Moreover, a generalization of Choquet theory developed by Mokobodoski and others \cite{Moko} yields that for every $\mu \in {\mathcal P}(X)$, we have 
\[
\sup \{\int_X f\ d\sigma;\, \mu \prec_{\mathcal A}\sigma \}=\int_X{\hat f}\, d\mu, 
\]
where   
\[
{\hat f}(x)=\inf \{ g(x); g\in -{\mathcal A},\,  g\geq f \, \hbox{on $X$}\}=\sup\{ \int_Xf d\sigma;\, \epsilon_x  \prec_{\mathcal A} \sigma\}.
\]
It follows that 
$
({\mathcal B}^b)_\mu^*(f)=\int_X{\hat f}\, d\mu$,  which means that ${\mathcal B}^b$ is a backward linear transfer whose Kantorovich operator is $T^-f={\hat f}$.
We can also define the {\it forward balayage transfer} as 
\begin{equation}
{\mathcal B}^f(\mu, \nu)=\left\{ \begin{array}{llll}
0 \quad &\hbox{if $\nu \prec_{\mathcal A} \mu$}\\
+\infty \quad &\hbox{\rm otherwise.}
\end{array} \right.
\end{equation} 
In this case, the forward Kantorovich operator is $T^+f=\check f$, where 
\[
{\check f}(x)=\sup \{ h(x); h\in {\mathcal A},\,  h \leq f \, \hbox{on $X$}\}=\inf\{ \int_Xf d\sigma;\, \epsilon_x  \prec_{\mathcal A} \sigma\}.
\]
\begin{itemize}
\item A typical example is when $X$ is a convex compact space in a locally convex topological vector space and ${\mathcal A}$ is the cone of continuous convex functions. In this case, $T^-f=\hat f$ (resp., $T^+f=\check f$)  is the concave envelope of $f$, and which was the context of the original Choquet theory. 

\item If $X$ is a bounded subset of a normed space $(E, \|\cdot\|)$, then ${\mathcal A}$ can be taken to be the cone of all norm-Lipschitz convex functions.  

\item If $X$ is an interval of the real line, then one can consider ${\mathcal A}$ to be the cone of increasing functions.

\item If $X$ is a pseudo-convex domain of $\C^n$, then one can take ${\mathcal A}$ to be the cone of Lipschitz plurisubharmonic functions (see \cite{GM}). In this case, if $\phi$ is a Lipschitz function, then the Lipschitz plurisubharmonic envelope of $\phi$, i.e., the largest Lipschitz PSH function below $\phi$ is given by the formula
\[
{\check \phi}(x)=\inf\{\int_0^{2\pi}\phi (P(e^{i\theta})\, \frac{d\theta}{2\pi}; P:\C\to X \,  \hbox{polynonial with}\,  P(0)=x\}.
\] 
Note that $\hat \phi=-\check \psi$, where $\psi=-\phi$. \\
\end{itemize}

 \noindent {\it Example 6: The Skorokhod transfer}
 
  Again, letting ${\mathcal S}$ be the class of --possibly randomized-- Brownian stopping times, and  
     define 
 \begin{equation}
{\mathcal SK}(\mu, \nu)=\left\{ \begin{array}{llll}
0 \quad &\hbox{if $B_0\sim \mu$ and $B_\tau \sim\nu$ for some $\tau \in {\mathcal S}$, }\\
+\infty \quad &\hbox{\rm otherwise.}
\end{array} \right.
\end{equation}  
The following is a  classical result of Skorokhod. See, for example \cite{G-K-L2} for a proof in higher dimension.
\begin{prop} Let ${\mathcal A}$ be the cone of Lipschitz subharmonic functions on a domain $\Omega$ in $\R^n$. Then, the following are equivalent for two probability measures $\mu$ and $\nu$ on $\Omega$.
\begin{enumerate}
\item $\mu \prec_{\mathcal A} \nu$ (i.e, $\mu$ and $\nu$ are in subharmonic order).
\item There exists a stopping time $\tau \in {\mathcal S}$ such that $B_0\sim \mu$ and $B_\tau \sim\nu$. 
\end{enumerate}
\end{prop}
This means that ${\mathcal SK}$ is a backward linear transfer with Kantorovich operator given by $T^-f=f_{**}$, which is the smallest Lipschitz superharmonic function above $f$. This can also be written as  $T^-f=J_f$, where $J_f(x)$ is a viscosity solution for the heat variational inequality,
 	\begin{equation}
		 \max\left\{f(x)-J(x), \Delta J(x)\right\}=0.
	\end{equation}
Another representation for $J_f$ is given 
  by  the following  dynamic programming principle,
	\begin{align} 
		J_f(x) := \sup_{\tau \in \mathcal S}\mathbb{E}^{x}\Big[f(B_\tau)\Big].
	\end{align}
	
 \subsection{Cost optimizing mass transports are backward and forward linear transfers}
 
The examples in this subsection correspond to cost minimizing transfers, where a cost $c(x,y)$ of moving two states is given. \\

   \noindent {\it Example 7: Monge-Kantorovich transfers} 
  
Any function $c\in C(X\times Y)$ determines a backward and forward linear transfer. This is Monge-Kantorovich theory of optimal transport. One associates the map ${\calT}_c$ on  ${\mathcal P}(X)\times  {\mathcal P}(Y)$ to be the optimal mass transport  between two probability measures $\mu$ on $X$ and $\nu$ on  $Y$, that is \begin{eqnarray}
{\mathcal T}_c(\mu, \nu):=\inf\big\{\int_{X\times Y} c(x, y)) \, d\pi; \pi\in \mK(\mu,\nu)\big\},
\end{eqnarray}
where $\mK(\mu,\nu)$ is the set of probability measures $\pi$ on $X\times Y$ whose marginal on $X$ (resp. on $Y$) is $\mu$ (resp., $\nu$) {\it (i.e., the transport plans)}. Monge-Kantorovich theory readily yields  that ${\mathcal T}_c$ is a linear transfer. Indeed, 
if we define the operators 
\begin{equation}
T ^+_cf(y)=\inf_{x\in X} \{c(x, y)+f(x)\} \quad {\rm and} \quad T ^-_cg(x)=\sup_{y\in Y} \{g(y)-c(x, y)\}, 
\end{equation}
for any $f\in C(X)$ (resp., $g\in C(Y)$), then Monge-Kantorovich duality yields that for any probability measures $\mu$ on $X$ and $\nu$ on  $Y$, we have 
\begin{eqnarray*}
{\mathcal T}_c(\mu, \nu)&=& 
\sup\big\{\int_{Y}T ^+_cf(y)\, d\nu(y)-\int_{X}f(x)\, d\mu(x);\,  f \in C(X)\big\}\\
&=&\sup\big\{\int_{Y}g(y)\, d\nu(y)-\int_{X}T _c^-g(x)\, d\mu(x);\,  g \in C(Y)\big\}.
\end{eqnarray*}
This means that the Legendre transform $({\mathcal T}_c)^*_\mu (g)=\int_{X}T _c^-g(x)\, d\mu(x)$ and $T_c^-$ is the corresponding backward Kantorovich operator. Similarly,  $({\mathcal T}_c)^*_\nu (f)=-\int_Y T ^+_c(-f)(y) \, d\nu(y)$ on $C(X)$ and $T ^+_c$ is the corresponding forward Kantorovich operator. See for example Villani \cite{V}. \\

\noindent {\it Example 7.1: The trivial Kantorovich transfer}
 
  Any pair of functions $c_1\in C(X)$, $c_2\in C(Y)$ defines trivially a linear transfer via
$$
{\cal T}(\mu, \nu)=\int_Yc_2\, d\nu -\int_X c_1\, d\mu.
$$
The Kantorovich operators are then
$T^+f=c_2+\inf (f-c_1)$  
and  $T^-g=c_1+\sup(g-c_2)$.\\
 
 \noindent {\it Example 7.2:  The Csisz\'ar-Kullback-Pinsker transfer} 
 
 This is simply the total variation distance between two probability measures $\nu$ and $\mu$ on $X$,  defined by
\begin{equation}
\|\nu-\mu\|_{TV}=\sup\{|\nu(A)-\mu(A)|; \hbox{$A$ measurable subset of $X$}\},
\end{equation}
with forward (resp., backward) Kantorovich operator given by
\begin{equation}
\hbox{$T^+f(y)=\min \{\inf\limits_{x\neq y}f(x)+1, f(y)\}$, while  $T^-g(x)=\max\{\sup\limits_{x\neq y}g(y)-1, g(x)\}$.}
\end{equation}
 It is actually a cost minimizing optimal transport, where the cost is given by the Hamming metric. \\
 
 \noindent {\it Example 7.3: The Kantorovich-Rubinstein transfer}
 
  If $d: X\times X\to \R$ is a lower semi-continuous metric on $X$, then 
 \begin{equation}
  {\cal T}(\mu,\nu) = \|\nu-\mu\|^*_\mathrm{Lip}:=\sup\left\{\int_X u \, d(\nu-\mu); u
\mathrm{\ measurable}, \|u\|_\mathrm{Lip}\leq 1\right\}
 \end{equation}
 is a linear transfer, where here    $
\|u\|_\mathrm{Lip}:=\sup_{x\not =y}\frac{|u(y)-u(x)|}{d(x,y)}$. The corresponding forward Kantorovich operator is then the Lipshitz regularization $T^+f(x)=\inf\{f(y)+d(y,x); y\in X\}$, while $T^-f(x)=\sup\{f(y)-d(x,y); y\in X\}$.  Note that $T^+\circ T^-f=T^-f$. \\

 \noindent {\it Example 7.4: The Brenier-Wasserstein distance} \cite{Br}
 
  If $c(x, y)=\langle x, y\rangle$ on $\R^d\times \R^d$, and  $\mu, \nu$ are two probability measures of compact support on $\R^d$, then 
 \begin{eqnarray*}
{\cal W}_2(\mu, \nu)=\inf\big\{\int_{\R^d\times \R^d} \langle x,y\rangle \, d\pi; \pi\in \mK(\mu,\nu)\big\}.
\end{eqnarray*}
Here, the Kantorovich operators are  
\begin{equation}T ^+f(x) 
=-f^*(-x) \quad \hbox{and \quad $T ^-g(y)=  (-g)^*(-y)$,}
\end{equation}
 where $f^*$ is the convex Legendre transform of $f$. \\

 \noindent {\it Example 7.5:  Optimal transport for a cost given by a generating function} (Bernard-Buffoni \cite{B-B1})
 
  This important example links the Kantorovich backward and forward operators with the forward and backward Hopf-Lax operators that solve first order Hamilton-Jacobi equations. Indeed, on a given compact manifold $M$, consider the cost:
\begin{equation}\label{BB}
c^L(y,x):=\inf\{\int_0^1L(t, \gamma(t), {\dot \gamma}(t))\, dt; \gamma\in C^1([0, 1), M);  \gamma(0)=y, \gamma(1)=x\},
\end{equation}
where $[0, 1]$ is a fixed time interval, and $L: TM \to \R\cup\{+\infty\}$ is a given Tonelli Lagrangian that is convex in the second variable of the tangent bundle $TM$. If now $\mu$ and $\nu$ are two probability measures on $M$, then 
 \begin{eqnarray*}
{\mathcal T}_{L}(\mu, \nu):=\inf\big\{\int_{M\times M} c^L(y, x) \, d\pi; \pi\in \mK(\mu,\nu)\big\}
 \end{eqnarray*}
is a linear transfer with forward Kantorovich operator given by 
 $T ^+_1f(x)=V_f(1, x)$, where
  $V_f(t, x)$ being the value functional
\begin{equation}\label{value.1}
V_f(t,x)=\inf\Big\{f(\gamma (0))+\int_0^tL(s,\gamma (s), {\dot \gamma}(s))\, ds; \gamma \in C^1([0, 1), M);   \gamma(t)=x\Big\}.
\end{equation}
Note that $V_f$ is --at least formally-- a solution for the Hamilton-Jacobi equation 
\begin{eqnarray}\label{HJ.0} 
\left\{ \begin{array}{lll}
\partial_tV+H(t, x, \nabla_xV)&=&0 \,\, {\rm on}\,\,  [0, 1]\times M,\\
\hfill V(0, x)&=&f(x). 
\end{array}  \right.
 \end{eqnarray}
 Similarly, the backward Kantorovich potential is given by  
 $T ^-_1g(y)=W_g(0, y),$ $W_g(t, y)$ being the value functional
 \begin{equation}\label{value.2}
W_g(t,y)=\sup\Big\{g(\gamma (1))-\int_t^1L(s,\gamma (s), {\dot \gamma}(s))\, ds; \gamma \in C^1([0, 1), M);   \gamma(t)=y\Big\},
\end{equation}
 which is a solution for the backward Hamilton-Jacobi equation 
\begin{eqnarray}\label{HJ.0} 
\left\{ \begin{array}{lll}
\partial_tW+H(t, x, \nabla_xW)&=&0 \,\, {\rm on}\,\,  [0, 1]\times M,\\
\hfill W(1, y)&=&g(y). 
\end{array}  \right.
 \end{eqnarray}
 
\subsection{One-sided linear transfers arising from constrained mass transports}
 
We now give examples of linear transfers, which do not fit in the framework of Monge-Kantorovich theory.  Cost minimizing mass transport with additional constraints give examples of one-directional linear transfers. We single out the following: \\

 \noindent {\it Example 8:  Martingale transports are backward linear transfers} 
  
  Martingale transports are $\cal C$-dilations when $\cal C$ is the cone of convex continuous functions on $\R^n$. In this case, $\mu\prec_C\nu$ is sometimes called the {\it Choquet order for convex functions}. Note that $x$ is the barycenter of a measure $\nu$ if and only if $\delta_x \prec_C\nu$, where $\delta_x$ is Dirac measure at $x$. \\
  So if $\mu, \nu$ are two probability measures, we then consider  $MT(\mu,\nu)$ to be the susbet of  ${\mathcal K}(\mu, \nu)$ consisting of {\it the martingale transport plans}, that is the set of probabilities $\pi$
on $\R^d \times \R^d$ with marginals $\mu$ and $\nu$, such that for $\mu$-almost $x\in\R^d$, 
the component $\pi_x$ of its disintegration $(\pi_x)_x$ with respect to $\mu$, i.e. $d\pi(x,y)=d\pi_x(y)d\mu(x)$, has its barycenter at $x$. As mentioned above, 
\begin{equation}
\hbox{$MT(\mu,\nu)\neq \emptyset$ if and only if $\mu \prec_{\cal C}\nu$.}
\end{equation}
 One can also use the probabilistic notation, which amounts to 
minimize 
$\E_{\rm P} \,c(X,Y)$
over all martingales $(X,Y)$ on a probability space $(\Omega, {\mathcal F}, P)$ into $\R^d \times \R^d$ (i.e. $E[Y|X]=X$) with laws $X \sim \mu$ and $Y \sim \nu$ (i.e., $P(X\in A)=\mu (A)$ and $P(Y\in A)=\nu (A)$ for all Borel set $A$ in $\R^d$). Note that in this case, the disintegration of $\pi$ can be written as the conditional probability   $\pi_x (A) =  \P {Y\in A|X=x}$.

If now $c:\R^d\times \R^d \to \R$ is a continuous cost function, then the corresponding martingale transport is a {\it backward linear transfer} in the following way:
\begin{equation}
{\mathcal T}_M(\mu, \nu)=\left\{ \begin{array}{llll}
\inf\{\int_{\R^d\times \R^d} c(x,y) \,d\pi(x,y); \pi\in MT(\mu,\nu)\} \quad &\hbox{if $\mu\prec_C\nu$}\\
+\infty \quad &\hbox{\rm if not.}
\end{array} \right.
\end{equation}
The backward Kantorovich operator is then given by
\[
\hbox{$T _M^-f(x)=\hat{f}_{c, x}(x)$, where $\hat{f}_{c, x}$ is the concave envelope of the function $f_{c, x}:y\to f(y)+c(x, y)$.}
\]
 See Henri-Labord\`ere \cite{HL} and Ghoussoub-Kim-Lim \cite{G-K-L1} for higher dimensions.\\

 \noindent {\it Example 9:  General stochastic transports are backward linear transfers} (Mikami-Thieulin \cite{M-T})
 
 Given a Lagrangian $L:[0,1]\times \R^d\times \R^d \rightarrow\mathbb{R}$, we define the following stochastic counterpart of the optimal transportation problem mentioned above.
  \begin{equation}
    {\mathcal T}_L(\mu,\nu):=\inf\left\{ \expect{\int_0^1 L(t,X(t),\beta_X(t,X(t)))\,dt}\middle\rvert X(0)\sim \mu,X(1)\sim\nu,X(\cdot)\in \mathcal{A}\right\}\\\label{eq:bstoc}
\end{equation}
Here $\mathcal{A}$ refers to the set of $\mathbb{R}^d$-valued continuous semimartingales $X(\cdot)$ such that there exists a measurable drift $\beta_X:[0,T]\times C([0,1])\rightarrow M^\ast$ where
\begin{itemize}
    \item $\omega\mapsto\beta_X(t,\omega)$ is $\mathcal{B}(C([0,t]))_+$-measurable for all $t$.
    \item $W(t):=X(t)-X(0)-\int_0^t \beta_X(s,X)\,ds$ is a $\sigma[X(s):s\in[0,t]]$-Brownian motion.
\end{itemize}
This stochastic transport does not fit in the standard optimal mass transport theory since it does not originate in optimization a cost between two deterministic states. However, under certain conditions on the Lagrangian, Mikami and Thieulin \cite{M-T} proved that the map $(\mu, \nu)\to   {\mathcal T}_L(\mu,\nu)$ is  jointly convex and weak$^*$-lower semi-continuous on the space of measures and that 
 \begin{equation}
    {\mathcal T}_L(\mu,\nu)=\sup\left\{\int_M f(x)\,d\nu-\int_M V_f(0,x)\,d\mu;\,  f\in\mathcal{C}_b^\infty\right\},
\end{equation}
where $V_f$ solves the Hamilton-Jacobi-Bellman equation
\begin{equation}\tag{HJB}
    \pderiv{V}{t}+\frac{1}{2}\Delta V(t,x)+H(t,x,\nabla V)=0, \quad V(1,x)=f(x).
\end{equation}
In other words,  ${\mathcal T}_L$ is a backward linear transfer with a Kantorovich operator being
 $T ^-_{L}f=V_f(0, \cdot), $
where $V_f(t, x)$ can be written as  
\begin{equation}
   V_f(t,x)=\sup_{X\in\mathcal{A}}\left\{\expcond{f(X(1))-\int_t^1 L(s,X(s),\beta_X(s,X))\,ds}{X(t)=x}\right\}.
\end{equation}

 \noindent {\it Example10:  Optimally stopped stochastic transport are backward linear transfers} (Ghoussoub-Kim-Palmer \cite{G-K-P1, G-K-P2})
 
 Consider the optimal stopping problem
 \begin{equation}
    {\mathcal T}_L(\mu,\nu)=\inf\left\{ \expect{\int_0^T L(t,X(t),\beta_X(t,X(t)))\,dt};X(0)\sim \mu, T \in {\cal S},  X(T)\sim\nu,X(\cdot)\in \mathcal{A}\right\},  
\end{equation}
where  ${\cal S}$ is the set of possibly randomized stopping times.
In this case, ${\mathcal T}_L$ is a backward linear transfer with Kantorovich potential given by 
$T ^-_{L}f={\hat V}_f(0, \cdot)$, where 
\begin{equation}
  {\hat V}_f(t,x)=\sup_{X\in {\mathcal A}}\sup_{T\in {\cal S}}\left\{\expcond{f(X(T))-\int_t^T L(s,X(s),\beta_X(s,X))\,ds}{X(t)=x}\right\},
\end{equation}
which is --at least formally-- a solution ${\hat V}_f(t, x)$ of  the quasi-variational Hamilton-Jacobi-Bellman inequality, 
	\begin{align} \label{eqn:HJB_0}
		\min\left\{\begin{array}{r} V_f(t,x)-f(x), \\ -\partial_t V_f(t,x)-{H}\big(t,x,\nabla V_f(t,x)\big)-\frac{1}{2}\Delta V_f(t,x)\end{array}\right\}=0.
	\end{align}
	
 \noindent {\it Example 11:  Optimal Skorokhod embeddings are backward linear transfers} (Ghoussoub-Kim-Palmer \cite{G-K-P3})
 
 Let $L(t, x)$ be a Lagrangian only depending on time and space and consider for any Radon probability measures $\mu$, $\nu$ with finite expectations,  
 the correalation
	\begin{align} \label{eqn:Skorokhod_cost}
		{\mathcal T}(\mu,\nu) :=  \inf
		\Big\{\mathbb{E}\Big[ \int_0^\tau L(t,B_t)dt\Big];\ \tau \in {\mathcal S}(\mu, \nu)
		\Big\},  
	\end{align}
where ${\mathcal S}(\mu, \nu)$ denotes the set of --possibly randomized-- stopping times with finite expectation such that $\nu$ is realized by the distribution of $B_\tau$ (i.e, $B_\tau \sim\nu$ in our notation), where $B_t$ is Brownian motion starting with $\mu$ as a source distribution, i.e., $B_0\sim \mu$. We shall assume ${\mathcal T}(\mu,\nu)=+\infty$ if ${\mathcal S}(\mu, \nu)=\emptyset$, which is the case if and only if $\mu$ and $\nu$ are not in subharmonic order. 
It has been proved in \cite{G-K-P3} that under suitable conditions, 
\begin{align} \label{eqn:dual_cost}
		{\mathcal T}(\mu,\nu) := \sup_{\psi}\Big\{\int_{\R^d}\psi(z)\nu(dz) - \int_{\R^d}J_\psi(0,y)\mu(dy);\ \psi\in C(\R^d)\Big\}, 
	\end{align}	
 where $J_\psi:\R^+\times \R^d\rightarrow \R$ is defined via the 
	dynamic programming principle
	\begin{align}\label{eqn:dynamic_programming}
		J_\psi(t,x) := \sup_{\tau \in \mathcal{R}^{t,x}}\Big\{\mathbb{E}^{t,x}\Big[\psi(B_\tau)-\int_t^\tau L(s,B_s)ds\Big]\Big\},
	\end{align}
	where the expectation superscripted with $t,x$ is with respect to the Brownian motions satisfying $B_t=x$, and the minimization is over all finite-expectation stopping times $\mathcal{R}^{t,x}$ on this restricted probability space such that $\tau \ge t$. In other words, $T^-\psi=J_\psi(0, \cdot)$ is a backward linear transfer. Note that $J_\psi(t, x)$ can be see as a ``variational solution" for 
the quasi-variational Hamilton-Jacobi-Bellman equation:
	\begin{align} 
		 \min\left\{\begin{array}{r} J(t,x) -\psi(x)\\ -\frac{\partial}{\partial t}J(t,x)-\frac{1}{2}\Delta J(t,x)+L(t,x)\end{array}\right\}=0.
	\end{align}

 \subsection{Weak optimal transports}

\noindent Other examples of linear transfers arise from the work of Marton, who extended the work of Talagrand.  \\

 \noindent {\it Example 12: Marton transports are backward linear transfers} (Marton  \cite{Ma1, Ma2})

These are transports of the following type: 
 \begin{equation}
 {\mathcal T}_{\gamma, d}(\mu, \nu)=\inf\left\{\int_X \gamma \left(\int_Y d(x,y)d\pi_x(y)\right)\, d\mu (x); \pi\in {\cal K}(\mu, \nu)\right\},
 \end{equation}
 where $\gamma$ is a convex function on $\R^+$ and $d:X\times Y \to \R$ is a lower semi-continuous functions. Marton's weak transfer correspond to $\gamma (t)=t^2$ and $d(x, y)=|x-y|$, which in probabilistic terms reduces to  
 \begin{equation}
 {\mathcal T}_2(\mu, \nu)=\inf\left\{\ensuremath{\mathbb{E}}[\ensuremath{\mathbb{E}}[|X-Y|\, |Y]^2]; X \sim \mu, Y\sim \nu \right\}. 
 \end{equation}
 This is a backward linear transfer with Kantorovich potential 
 \[
 T ^-f(x)=\sup\left\{ \int_Yf(y) d\sigma (y) - \gamma \left(\int_Y d(x, y)\, d\sigma (y)\right); \ \sigma \in {\cal P}(Y)\right\}.
 \]
  
 \noindent {\it Example 12.1:  A barycentric cost function} (Gozlan et al. \cite{Go4})
  
Consider the (weak) transport
  \begin{equation}
  {\cal T}(\mu, \nu)=\inf\left\{\int_X \|x-\int_Y y d\pi_x(y)\|\, d\mu (x); \pi\in {\cal K}(\mu, \nu)\right\}.
 \end{equation}
Again, this is a backward linear transfer, with Kantorovich potential 
 \[
 T ^-f(x)=\sup\{ f_{**} (y) -\|y-x\|; y\in \R^n\}.
 \]
where $ f_{**}$ is the concave envelope of $f$, i.e., the smallest concave usc function above $f$. \\
   
  \noindent {\it Example 12.2: Schr\"odinger bridge} (Gentil-Leonard-Ripani  \cite{GLR})
 
Fix some reference non-negative measure $R$ on path space $\Omega=C([0,1], \R^n)$, and let $(X_t)_t$ be a random process on $\R^n$ whose law is $R$. Denote by $R_{01}$ the joint law of the initial position $X_0$ and the final position $X_1$, that is $R_{01}=(X_0, X_1)_\#R$. For probability measures $\mu$ and $\nu$ on $\R^n$, the  maximum entropy formulation of the Schr\"odinger bridge problem between $\mu$ and $\nu$ is defined as 
 \begin{equation}
 {\cal S}_R(\mu, \nu)=\inf\{\int_{\R^n \times \R^n}\log(\frac{d\pi}{dR_{01}})\, d\pi; \pi \in {\cal K}(\mu, \nu\}.
 \end{equation} 
For example (See \cite{GLR}), assume $R$ is the reversible Kolmogorov continuous Markov process associated with the generator 
 $\frac{1}{2}(\Delta -\nabla V\cdot \nabla)$ and the initial mesure $m=e^{-V(x)}dx$ for some function $V$. Then, under appropriate conditions on $V$ (e.g., if $V$ is uniformly convex), then 
 \[
 {\cal T}(\mu, \nu)={\cal S}_R(\mu, \nu)-\frac{1}{2}\int_{\R^n}\log(\frac{d\mu}{dm})\, d\mu-\frac{1}{2}\int_{\R^n}\log(\frac{d\nu}{dm})\, d\nu
 \] 
 is a forward linear transfer with Kantorovich operator 
 \[
 T ^+f(x)=\log E_{R^x}e^{f(X_1)}=\log S_1(e^f)(x),
 \]
 where $(S_t)$ is the  semi-group associated to $R$. 
  It is worth noting that $ {\cal T}$ is symmetric, that is ${\cal T}(\mu, \nu)= {\cal T}(\nu, \mu)$, which means that it is also a backward linear transfer.  Note that when $V=0$, the process is Brownian motion with Lebesgue measure as its initial reversing measure, while when $V(x)=\frac{|x|^2}{2}$, $R$ is the path measure associated with the Ornstein-Uhlenbeck process withe the Gaussian as its initial reversing measure. 

 \section{A representation of linear transfers as generalized optimal mass transports}
 
  We now consider whether any transfer  ${\mathcal T}$ on $X\times Y$ arises from a cost minimizing mass transport. Note first that transfers need not be  defined on Dirac measures, a prevalent situation in stochastic transport problems. Moreover, even if the set of Dirac measures $\{(\delta_x, \delta_y); (x, y)\in X\times Y\} \subset D({\calT})$, and we can then define a cost function as $c(x,y)=\calT (\delta_x, \delta_y)$, and its associated optimal mass transport  ${\mathcal T}_c(\mu, \nu)$, we then only have  
 \begin{equation}\label{compare}
{\mathcal T}_c(\mu, \nu) \geq {\mathcal T}(\mu, \nu).
\end{equation}
Indeed, for every $x\in X$, we have 
\[
\hbox{${T ^-}g (x)={\calT}_{\delta_x}^*(g)=\sup\{\int_Y g d\nu -{\mathcal T}(\delta_x, \nu); \, \nu\in {\mathcal P}(Y) \}\geq \sup\{g(y) -c(x,y); y\in Y\}=T _c^-g(x)
$, }
\]
hence, 
\begin{eqnarray*}
{\mathcal T}(\mu, \nu) 
=\sup\{\int_Y gd\nu  -\int_X {T ^-}g \, d\mu;\, g\in C(Y)\}\leq 
\sup\{\int_Y gd\nu  -\int_X {T ^-_c}g \, d\mu;\, g\in C(Y)\}= {\mathcal T}_c(\mu, \nu).
\end{eqnarray*}
Moreover, the inequality (\ref{compare}) is often strict.

 Motivated by the work of Marton and others, Gozlan et al. \cite{Go4} introduced the notion of {\it weak transport}. It consists of  considering cost minimizing transport plans, where cost functions between two points are replaced by  {\it generalized costs} $c$ on $X\times {\mathcal P}(Y)$, 
 where $\sigma \to c(x, \sigma)$ is convex and lower semi-continuous. As the following proposition shows. this notion turns out to be equivalent to the notion of backward linear transfer, at least in the case where Dirac measures belong to the first partial effective domain of the map $\T$, that is when $\{\delta_x; x\in X\}\subset D_1 ({\mathcal T})$.

\begin{prop} Let  ${\mathcal T}: {\mathcal P}(X)\times {\mathcal P}(Y)\to \R \cup\{+\infty\}$ be a functional  
such that $\{\delta_x; x\in X\}\subset D_1 ({\mathcal T})$. Then, ${\cal T}$ is a backward linear transfer  if and only if there exists a lower semi-continuous function $c: X\times {\mathcal P}(Y)\to \R\cup \{+\infty\}$ with $\sigma \to c(x, \sigma)$ convex on ${\mathcal P}(Y)$ for each $x\in X$ such that for every $\mu \in {\mathcal P}(X)$, and $\nu \in {\mathcal P}(Y)$, we have
  \begin{equation}
 {\mathcal T}(\mu, \nu)=\inf_\pi\{\int_Xc(x, \pi_x)\, d\mu(x); \pi \in {\mathcal K}(\mu, \nu)\}.
 \end{equation}
 The corresponding backward Kantorovich operator is given for every $g\in C(Y)$ by
 \begin{equation}\label{Kminus}
 T ^-g(x)=\sup\{\int_Yg(y)\, d\sigma (y)-{\mathcal T}(x, \sigma); \sigma \in {\mathcal P}(Y)\}.
 \end{equation} 
  \end{prop}
Note that we have identified here any $\pi \in {\mathcal K}(\mu, \nu)$ with its disintegration that gives a probability kernel $\pi:X\to {\mathcal P}(X)$ such that 
 $\nu (A)=\int_X \pi_x (A) d\mu (x)$. \\
 
\noindent{\bf Proof:} Consider first a lower semi-continuous function $c: X\times {\mathcal P}(Y)\to \R\cup \{+\infty\}$ with $\sigma \to c(x, \sigma)$ convex on ${\mathcal P}(Y)$ for each $x\in X$, and let
\[
{\mathcal T}_c(\mu, \nu):=\inf_\pi\int_Xc(x, \pi_x)\, d\mu(x); \pi \in {\mathcal K}(\mu, \nu)\}.
\]
We first prove that ${\mathcal T}_c$ is a backward linear transfer  with a Kantorovich operator given by 
\begin{equation}\label{K}
 T_c ^-g(x)=\sup\{\int_Yg(y)\, d\sigma (y)-c(x, \sigma); \sigma \in {\mathcal P}(Y)\}.
 \end{equation} 
  This will then imply that if ${\mathcal T}$ is any backward linear transfer with Kantorovich operator $T^-$, and if $c(x, \sigma)={\mathcal T}(\delta_x, \sigma)$, then $
T ^-g(x)= T_c^-g(x)$ and therefore
 $
  {\mathcal T}(\mu, \nu)= {\mathcal T}_c(\mu, \nu).
  $\\
  First, it is easy to show that $ {\mathcal T}_c$ is a convex lower semi-continuous function on ${\mathcal P}(X) \times {\mathcal P}(Y)$. Consider now the Legendre transform of $({\mathcal T}_c)_\mu$, that is 
 \begin{eqnarray*}
 ({\mathcal T}_c)^*_\mu (g)&=&\sup\{\int_Yg\, d\nu -  {\mathcal T}_c(\mu, \nu); \nu\in {\mathcal P}(Y)\}\\
 &=&\sup\{\int_Yg(y)\, d\nu (y) -  \int_Xc(x, \pi_x)\, d\mu(x); \nu\in {\mathcal P}(Y), \pi \in {\mathcal K}(\mu, \nu)\}\\
&=&\sup\{\int_X\int_Yg(y) d\pi_x(y)\, d\mu (x) -  \int_Xc(x, \pi_x)\, d\mu(x);  
\pi \in {\mathcal K}(\mu, \nu)\}\\
&\leq &\sup\{\int_X\int_Yg(y) d\sigma(y)\, d\mu (x) -  \int_Xc(x, \sigma)\, d\mu(x); \sigma \in {\mathcal P}(Y)\}\\
&\leq &\int_X\{\sup\limits_{\sigma \in {\mathcal P}(Y)}\{\int_Yg(y) d\sigma(y) - c(x, \sigma)\}\, d\mu(x)\}\\
&=&\int_XT_c^-g(x) d \mu (x).
 \end{eqnarray*}
On the other hand, use your favorite selection theorem to find a measurable selection $x\to {\bar \pi}_x$ from $X$ to ${\mathcal P}(Y)$ such that
\[
T_c^-g(x)=\int_Yg(y) d {\bar \pi}_x(y) - c(x, \pi_x) \quad \hbox{for every $x\in X$.}
\]
It follows that 
 \begin{eqnarray*}
({\mathcal T}_c)^*_\mu (g)&=&\sup\{\int_Yg\, d\nu -  {\mathcal T}_c(\mu, \nu); \nu\in {\mathcal P}(Y)\}\\
&\geq & \int_X\{\int_Yg(y) d\pi_x(y) - c(x, \pi_x)\}\, d\mu (x)\\
&=& \int_XT^-_cg(x) d \mu (x), 
 \end{eqnarray*}
hence   $({\mathcal T}_c)^*_\mu (g)=\int_XT^-_cg(x) d \mu (x)$ and  $T^-=T_c^-$.\\
Conversely, if ${\mathcal T}$ is a backward linear transfer  with $T ^-$ as a Kantorovich operator, then by setting $c(\delta_x, \sigma)={\mathcal T}(\delta, \sigma)$, we have 
$T ^-=T _c^-g$ and we are done.

  \section{Operations on linear mass transfers}
 
 Denote by ${\cal LT}_- (X \times Y)$ (resp., ${\cal LT}_+ (X \times Y)$) the class of backward (resp., forward) linear transfers on $X\times Y$. The following proposition is an immediate consequence of the properties of the Legendre transform.
  \begin{prop} The class ${\cal LT}_- (X \times Y)$ (resp., ${\cal LT}_+f (X \times Y)$) is a convex cone in the space of convex weak$^*$ lower continuous functions on ${\mathcal P}(X)\times {\mathcal P}(Y)$. 
 
 \begin{enumerate}
 
 \item {\bf (Scalar multiplication)} If $a\in \R^+\setminus \{0\}$ and ${\mathcal T}$ is a backward linear  
transfer with Kantorovich operator $T ^-$, 
then the transfer $(a{\mathcal T})$ defined by
$(a{\mathcal T})(\mu, \nu)=a{\mathcal T}(\mu, \nu)$ is also a backward linear 
transfer with Kantorovich operator on $C(Y)$ defined by,
\begin{equation}
T _a^-(f)=aT ^-(\frac{f}{a}).
\end{equation}

\item {\bf (Addition)} If ${\mathcal T}_1$ and ${\mathcal T}_2$ are backward linear transfers  on $X\times Y$ with Kantorovich operator $T_1^-$, $T_2^-$ respectively, and such that $X\subset D({\mathcal T}_1)\cap D({\mathcal T}_2)$, then $ ({\mathcal T}_1 +  {\mathcal T}_2) (\mu, \nu):={\mathcal T}_1 (\mu, \nu)+  {\mathcal T}_2 (\mu, \nu)$
 is a backward linear transfer  on $X\times Y$, with Kantorovich operator given on $C(Y)$ by 
\begin{equation}
T^-f(x)=\inf\{T_1^-g(x)+T_2^-(f-g)(x); g\in C(Y)\}.
\end{equation}
\end{enumerate}
 Similar statements hold for ${\cal LT}_+f (X \times Y)$.
 \end{prop}

 \begin{defn} Consider the following operations on transfers. 

 \begin{enumerate}
 
\item  ({\bf Inf-convolution}) Let $X_1, X_2, X_3$ be 3 spaces, and suppose ${\mathcal T}_1$ (resp., ${\mathcal T}_2$) are functionals on  ${\mathcal P}( X_1)\times {\mathcal P}(X_2)$ (resp., ${\mathcal P}( X_2)\times {\mathcal P}(X_3)$). The {\it convolution} of ${\mathcal T}_1$ and ${\mathcal T}_2$ is the functional on ${\mathcal P}( X_1)\times {\mathcal P}(X_3)$ given by
 \begin{equation}
{\mathcal T}(\mu, \nu):={\mathcal T}_{1}\star{\mathcal T}_{2}=
\inf\{{\mathcal T}_{1}(\mu, \sigma) + {\mathcal T}_{2}(\sigma, \nu);\, \sigma \in {\mathcal P}(X_2)\}.
\end{equation}

\item ({\bf Tensor product}) If ${\mathcal T}_1$ (resp., ${\mathcal T}_2$) are functionals  on ${\mathcal P}( X_1)\times {\mathcal P}(Y_1)$ (resp., ${\mathcal P}( X_2)\times {\mathcal P}(Y_2)$) such that $X_1\subset D({\mathcal T}_1)$ and $X_2\subset D({\mathcal T}_2)$, then the tensor product of  ${\mathcal T}_1$ and ${\mathcal T}_2$ is the functional on ${\mathcal P}(X_1\times X_2)\times {\mathcal P}(Y_1\times Y_2)$ defined by:
\begin{eqnarray*}
 {\mathcal T}_1 \tens  {\mathcal T}_2 (\mu, \nu)=\inf\left\{ \int_{X_1\times X_2}\big({\mathcal T}_1 (x_1, \pi_{x_1,x_2})+{\mathcal T}_2 (x_2, \pi_{x_1,x_2})\big)\, d \mu (x_1, x_2); \pi \in {\mathcal K}(\mu, \nu)\right\}.
\end{eqnarray*}
\end{enumerate}
\end{defn}
The following easy proposition is important to what follows.

\begin{prop} \label{inf.tens} The following stability properties hold for the class of backward linear transfers.
\begin{enumerate}

\item If ${\mathcal T}_1$ (resp., ${\mathcal T}_2$) is a backward linear transfer  on  $X_1\times  X_2$ (resp., on $X_2\times X_3$) with Kantorovich operator $T _1^-$ (resp., $T _2^-$) , then ${\mathcal T}_{1}\star{\mathcal T}_{2}$ is also a backward linear transfer  on  $ X_1\times X_3$ with Kantorovich operator equal to $T _1^-\circ T _2^-$. 

\item If ${\mathcal T}_1$ (resp., ${\mathcal T}_2$) is a backward linear transfer  on $X_1\times Y_1$ (resp., $X_2\times Y_2$) such that $X_1\subset D({\mathcal T}_1)$ and $X_2\subset D({\mathcal T}_2)$, then $ {\mathcal T}_1 \tens  {\mathcal T}_2$
  is a backward linear transfer  on $(X_1\times X_2)\times (Y_1\times Y_2)$, with   Kantorovich operator given by 
\begin{equation}
T ^-g(x_1, x_2)=\sup\{\int_{Y_1\times Y_2} f(y_1, y_2) d\sigma(y_1, y_2) - {\mathcal T}_1(x_1, \sigma_1)-{\mathcal T}_2(x_2, \sigma_2);\,  \sigma \in {\mathcal K}(\sigma_1, \sigma_2)  \}.
 \end{equation} 
 Moreover, \begin{equation}\label{tens}
    {\mathcal T}_1 \tens  {\mathcal T}_2(\mu,\nu_1\otimes\nu_2)
 \leq
 \mathcal{T}_{1}(\mu_1,\nu_1)+\int_{X_1}\mathcal{T}_{2}(\mu_2^{x_1},\nu_2)\,d\mu_1(x_1), 
\end{equation}
where  $d\mu(x_1,x_2)=d\mu_1(x_1)d\mu_2^{x_1}(x_2).$
 
\end{enumerate}

\end{prop}
Note that a similar statement holds for forward linear transfers, modulo order reversals. For example,   if ${\mathcal T}_1$ and ${\mathcal T}_2$) are forward linear transfer, then  ${\mathcal T}_{1}\star{\mathcal T}_{2}$ is a  forward linear transferon  $ X_1\times X_3$ with Kantorovich operator equal to $T _2^+\circ T _1^+$.  

 \noindent{\bf Proof:}  For 1), we note first that if ${\mathcal T}_{1}$ (resp., ${\mathcal T}_{2}$) is jointly convex and weak$^*$-lower semi-continuous on ${\mathcal P}(X_1)\times {\mathcal P}(X_2)$ (resp., ${\mathcal P}(X_2)\times {\mathcal P}(X_3)$),  
 then both  $({\mathcal T}_{1}\star{\mathcal T}_{2})_\nu: \mu \to ({\mathcal T}_{1}\star{\mathcal T}_{2})(\mu, \nu)$ and  $({\mathcal T}_{1}\star{\mathcal T}_{2})_\mu: \nu \to ({\mathcal T}_{1}\star{\mathcal T}_{2})(\mu. \nu)$ are convex and weak$^*$-lower semi-continuous. We now calculate their Legendre transform. For $g\in C(X_3)$,
 \begin{eqnarray*}
 ({\mathcal T}_{1}\star{\mathcal T}_{2})_\mu^*(g)&=&\sup\limits_{\nu \in {\mathcal P}(X_3)}\sup\limits_{\sigma \in  {\mathcal P}(X_2)} 
 \left\{\int_{X_3} g\, d\nu -{\mathcal T}_{1}(\mu, \sigma) - {\mathcal T}_{2}(\sigma, \nu)\right\}\\
&=&  
\sup\limits_{\sigma \in  {\mathcal P}(X_2)} 
 \left\{({\mathcal T}_{2})_\sigma^* (g)-{\mathcal T}_{1}(\mu, \sigma) \right\}\\
&=& \sup\limits_{\sigma \in  {\mathcal P}(X_2)} 
 \left\{\int_{X_2} T _2^-(g)\, d\sigma-{\mathcal T}_{1}(\mu, \sigma) \right\}\\
 &=&({\mathcal T}_{1})_\mu^* (T _2^-(g))\\
 &=& \int_{X_1} T _1^-\circ T _2^-g\, d\mu.
 \end{eqnarray*} 
  In other words, 
$
 {\mathcal T}_{1}\star{\mathcal T}_{2}(\mu, \nu)=
\sup\big\{\int_{X_3}g(x)\, d\nu(x)-\int_{X_1} T _1^-\circ T _2^-g\, d\mu;\,  f\in C(X_3) \big\}
$.\\

 2) follows immediately from the last section since we are defining the tensor product as a generalized cost minimizing transport, where the cost ion $X_1\times X_2 \times {\mathcal P}(Y_1\times Y_2)$ is simply,
 \[
 {\mathcal T}((x_1, x_2),  \pi)= {\mathcal T}_1(x_1, \pi_1)+{\mathcal T}_2(x_1, \pi_2), 
 \]
where $\pi_1, \pi_2$ are the marginals of $\pi$ on $Y_1$ and $Y_2$ respectively. ${\mathcal T}_1 \tens  {\mathcal T}_2$ is clearly its corresponding backward transfer  with $T ^-$ being its Kantorovich operator. 

More notationally cumbersome but straightforward  is how to write the Kantorovich operators of the tensor product $T ^-g(x_1, x_2)$ in terms of $T_1^-$ and $T_2^-$, in order to establish (\ref{tens}). \\

 \noindent {\it Example 13:  Stochastic ballistic transfer} (Barton-Ghoussoub \cite{B-G})

Consider the stochastic ballistic transportation problem defined as:
\begin{equation}
    \underline{\mathcal B}(\mu,\nu):=\inf\left\{\expect{\langle V,X(0)\rangle +\int_0^T L(t,X,\beta_X(t,X))\,dt}\middle\rvert V\sim\mu, X(\cdot)\in \mathcal{A},X(T)\sim \nu\right\},
\end{equation}
where we are using the notation of Example 9.  Note that this a convolution of the Brenier-Wasserstein transfer of Example 7.3 with the general stochastic transfer of Example 9. Under suitable conditions on $L$, one gets that 
\begin{equation}
    \underline{\mathcal B}(\mu,\nu)=\sup \left\{\int g\, d\nu-\int \widetilde{\psi_g}\,d\mu; g\in C_b\right\},
\end{equation}
where $\widetilde{h}$ is the concave legendre transform of $-h$ and $\psi_{g}$ is the solution to the Hamilton-Jacobi-Bellman equation
\begin{align}\tag{HJB}
    \pderiv{\psi}{t}+\frac{1}{2}\Delta\psi(t,x)+H(t,x,\nabla\psi)=0,\quad  \psi(1,x)=g(x).
\end{align}
In other words, $ \underline{\mathcal B}$ is a backward linear transform with Kantorovich operator $T^-g=\widetilde{\psi_g}$. 
   
 \begin{rmk}{\rm (Lifting convolutions to Wasserstein space)} Let $X_0, X_1,...., X_n$ be compact spaces, and suppose for each $i=1,..., n$, we have a cost function $c_i: X_{i-1}\times X_i$, we consider the following cost function on $X_0\times X_n$, defined by
\[
c(y, x)=\inf\left\{c_1(y, x_1) +c_2(x_1, x_2) ....+c_n(x_{n-1}, x); \, x_1\in X_1, x_2\in X_2,..., x_{n-1}\in X_{n-1}\right\}. 
\]
 Let $\mu$ (resp., $\nu$ be probability measures on $X_0$ (resp., $X_n)$.  
\begin{enumerate}
\item The following then holds on Wasserstein space:
\begin{equation}
{\mathcal T}_c(\mu, \nu)=\inf\{{\mathcal T}_{c_1}(\mu, \nu_1) + {\mathcal T}_{c_2}(\nu_1, \nu_2) ...+{\mathcal T}_{c_n}(\nu_{n-1}, \nu);\, \nu_i \in {\mathcal P}(X_i), i=1,..., n-1\},
\end{equation}
and the infimum is attained at $\bar \nu_1, \bar \nu_2,..., \bar \nu_{n-1}$.
\item The following duality formula holds:
 \begin{eqnarray*}
{\mathcal T}_c(\mu, \nu)&=& 
\sup\big\{\int_{X_n}T^+_{c_n}\circ T^+_{c_{n-1}}\circ...T^+_{c_1}f(x)\, d\nu(x)-\int_{X_0}f(y)\, d\mu(y);\,  f\in C(X_0) \big\}\\
&=&\sup\big\{\int_{X_n}g(x)\, d\nu(x)-\int_{X_0}T ^-_{c_n}\circ T^-_{c_{n-1}}...\circ T^-_{c_1}g(x);\,  g\in C(X_n) \big\}.
\end{eqnarray*}

 \end{enumerate}

\end{rmk}

We note a few more elementary properties of linear transfers. 
\begin{prop} Let ${\mathcal T}: {\mathcal P}(X)\times {\mathcal P}(Y) \to \R\cup \{+\infty\}$ be a proper (jointly) convex and weak$^*$ lower semi-continuous on ${\mathcal M}(X)\times {\mathcal M}(Y)$.
 If  ${\mathcal T}$ is both a forward and backward linear transfer, 
and if  $\{(\delta_x, \delta_y); (x, y)\in X\times Y\} \subset D({\calT})$, then for any $g\in C(Y)$ such that $T^+g\in C(X)$ (resp.,$f\in C(X)$ such that $T^-f\in C(Y)$)  
\begin{equation}
\hbox{$T ^+\circ T ^-g(y) \geq g(y)$ for $y\in Y,$ \quad  \quad $T ^-\circ T ^+f(x) \leq f(x)$ for $x\in X$,}
\end{equation}
and
 \begin{equation} 
 \hbox{$T ^-\circ T ^+\circ T ^-g=  T ^-g$ and $T ^+\circ T ^-\circ T ^+g = T ^+g.$}
 \end{equation} 
\end{prop}
\noindent{\bf Proof:}  
Write for  $\nu \in {\mathcal P}(Y)$, 
\begin{eqnarray*}
\int_YT ^+\circ T ^-g\, d\nu&=&-T_\nu^*(-T ^-g)\\  
&=& -\sup\{-\int_XT_{\delta_x}^*(g)\, d\mu (x)-{\calT}(\mu, \nu); \mu \in {\mathcal P}(X)\}\\
&=& \inf \{\int_XT_{\delta_x}^*(g)\, d\mu (x)+{\calT}(\mu, \nu); \mu \in {\mathcal P}(X)\}\\
&\geq & \inf \{\int_XT_{\delta_x}^*(g)\, d\mu (x)+\int_Yg\, d\nu - \int_XT_{\delta_x}^*(g)\, d\mu; \mu \in {\mathcal P}(X)\}\\
&=&\int_Yg\, d\nu.
\end{eqnarray*}
 The last item follows from the above and the positivity property  of the Kantorovich operators.
\begin{rmk} By analogy with the case of cost optimizing mass transports, and assuming that a transfer is both forward and backward, we can say that a function $f$ in $C(X)$ is ${\cal T}$-concave if it is of the form $f=T^-\circ T^+g$ for some $g\in C(X)$.  It follows from the last proposition that 
 \begin{equation}
{\mathcal T}(\mu, \nu)=\sup\big\{\int_{Y}T^+f(y)\, d\nu(y)-\int_{X}f(x)\, d\mu(x);\,  f \in C(X),\hbox{ $f$ is ${\cal T}$-concave} \big\}.
\end{equation}

Similarly, we can say that a function $g$ in $C(Y)$ is ${\cal T}$-convex if it is of the form $g=T^+\circ T^-f$ for some $f\in C(Y)$.  
In this case, 
 \begin{equation}
{\mathcal T}(\mu, \nu)=\sup\big\{\int_{Y}g(y)\, d\nu(y)-\int_{X}T^-g\, d\mu(x);\,  g \in C(Y),\hbox{ $g$ is ${\cal T}$-convex} \big\}.
\end{equation}
\end{rmk}

  \subsection*{Duality for projections on subsets of Wasserstein space}
 
Let $\T$ be a linear transfer on $X\times Y$ and $K$ a closed convex set of probability measures on $Y$. We consider the following minimization problem 
\begin{equation}
 \inf\{ {\cal T}(\mu, \sigma); \sigma \in K\},
\end{equation}
which amounts to finding ``the projection" of $\mu$ on $K$, when the ``distance" is given by the transfer $\T$. In some cases, the set $K:={\cal C}(\nu)$ is a convex compact subset of $\cal P (Y)$ that depends on a probability measure $\nu$ in such a way that the following map
\begin{equation*}
{\cal S}(\sigma, \nu)=\left\{ \begin{array}{llll}
0 \quad &\hbox{if $\sigma \in {\cal C}(\nu)$}\\
+\infty \quad &\hbox{\rm otherwise.}
\end{array} \right.
\end{equation*}
is a backward transfer on $Y\times Y$. It then follows that 
 \[
 \inf\{ {\cal T}(\mu, \sigma); \sigma \in {\cal C}(\nu)\}= \inf\{ {\cal T}(\mu, \sigma)+ {\cal S}(\sigma, \nu); \sigma \in \cal P (X)\}={\cal T}\star {\cal S} (\mu, \nu).
 \]
If now $T^-$ (resp., $S^-)$ is the backward Kantorovich operator for ${\cal T}$ (resp., ${\cal S}$), then by Proposition \ref{con.con}, the Kantorovich operator for ${\cal T}\star {\cal S}$ is $T^-\circ S^-$, that is
 \begin{equation}
 \inf\{ {\cal T}(\mu, \sigma); \sigma \in {\cal C}(\nu)\}=\sup\{\int_Yg\, d\nu-\int_{X}T^-\circ S^-g\, d\mu;  g\in C(Y)\}.
\end{equation}
Here is an example motivated by a recent result of Gozlan-Juillet \cite{GoJu}. \\

  \noindent {\it Example 14: Projection on the set of balay\'ees of a given measure} 

Consider the problem 
\begin{equation}
{\cal P}(\mu, \sigma)=\inf\{ {\cal T}_c(\mu, \sigma); \sigma\prec_C \nu\}, 
\end{equation}
 where  ${\cal T}_c$ is the optimal mass transport associated to a cost $c(x, y)$ on $X\times Y$, and $\prec_C$ is the convex order on a convex compact set $Y$. Then, 
 \[
  {\cal P}(\mu, \nu)= {\cal T}_c\star {\mathcal B} (\mu, \nu)
\]
where ${\cal B}^b$ is the backward Balayage transfer. It follows that ${\cal P}$ is a backward linear transfer  
with Kantorovich operator being the composition of those for ${\cal T}_c$ and ${\cal B}^b$, that is
 \[
 T ^-f(x)=\sup\{ {\hat f} (y) - c(x,y); y\in Y\}.
 \]
where $\hat f$ is the concave envelope of $f$ on $Y$. 
We note that this is the same Kantorovich operator as for 
 the (weak) barycentric transport
  \begin{equation*}
  {\cal T}^c_B(\mu, \nu)=\inf\left\{\int_X c(x, \int_Y y d\pi_x(y))\, d\mu (x); \pi\in {\cal K}(\mu, \nu)\right\},
 \end{equation*}
 at least in the case where $c(x, y)=\alpha (x-y)$ for some convex lower semi-continuous $\alpha :\R \to \R^+$ (See \cite{Go4}). 
This then yields that ${\cal P}(\mu, \nu)= {\cal T}(\mu, \nu)$. 

Note that the martingale transport of Example 8 can be written as 
\[
{\cal T}^c_M={\cal T}_c+ {\cal B}^b,
\]
while the  (weak) barycentric transport
\[
  {\cal T}^c_{B}={\cal T}_c\star {\cal B}^b.
\]
Similar manipulations can be done when the balayage is given by the cones of subharmonic or plurisubharmonic functions. 

 \section{Examples of convex and entropic transfers}
 
  We now give a few examples of convex and entropic transfers, which are not necessarily linear transfers. First, recall that the  {\it increasing Legendre transform} (resp., {\it decreasing Legendre transform)} of a function $\alpha: \R^+\to \R$ (resp., $\beta: \R^+ \setminus \{0\} \to \R$) is defined as
    \begin{equation}
  \alpha^{\oplus}(t)=\sup\{ts-\alpha(s); s\geq 0\} \quad\hbox{resp.,  $\beta^{\ominus}(t)=\sup\{-ts-\beta(s); s > 0\}$ }
   \end{equation}
By extending $\alpha$ to the whole real line by setting $\alpha (t)=+\infty$ if $t<0$, and using the standard Legendre transform, one can easily show that $\alpha$ is convex increasing on $\R^+$ if and only if  $\alpha^{\oplus}$ is convex and increasing on $\R^+$. We then have the following reciprocal formula
  \begin{equation}
\alpha (t)=\sup\{ts-\alpha^{\oplus} (s); s\geq 0\}.
  \end{equation}
Similarly, if $\beta$ is convex decreasing on $\R^+\setminus \{0\}$, we have
 \begin{equation}
\beta (t)=\sup\{-ts-\beta^{\ominus} (s); s\geq 0\}.
  \end{equation}

  \begin{prop} Let $\alpha: \R^+\to \R$ (resp., $\beta: \R^+ \setminus \{0\} \to \R$) be a convex (resp., concave) increasing functions.
    \begin{enumerate}
   \item If ${\cal T}$ is a linear backward (resp., forward) transfer with convex Kantorovich operator $T^-$ (resp., concave Kantorovich operator $T^+$), then $\alpha ({\cal T})$ is a backward convex transfer (resp., forward convex transfer) with Kantorovich operators $(T_s^-)_{s\geq 0}$ (resp.,$(T_s^+)_{s\geq 0}$, where
\begin{equation}
T_s^-f=sT^-(\frac{f}{s})-\alpha^\oplus(s)\quad \hbox{(resp., $T_s^+f=sT^+(\frac{f}{s})-\alpha^\oplus(s)$.}
\end{equation}
In particular, for any $p\geq 1$, $\T^p$ is a convex forward (resp., backward) transfer.
 
   \item If ${\mathcal E}$ is a $\beta$-backward transfer     with Kantorovich operator $E ^-$, then it 
is a backward convex transfer with Kantorovich operators $(T^-_s)_{s>0}$ given by
\begin{equation}
T_s^-f=sT^-f+(-\beta)^\ominus(s).
\end{equation}

\item Similarly, if ${\mathcal E}$ is an $\alpha$-forward transfer     with Kantorovich operator $E ^+$, then it 
is a forward convex transfer with Kantorovich operators $(T^+_s)_{s>0}$ given by
\begin{equation}
T_s^+f=sT^+f-\alpha^\oplus(s).
\end{equation}
   \end{enumerate}
   \end{prop} 
\noindent{\bf Proof:}
For 1) it suffices to write 
\begin{eqnarray*}
\alpha ({\mathcal T}(\mu, \nu))
&=& 
\sup\big\{s\int_{Y}{T^+}f\, d\nu-s\int_{X}f\, d\mu -\alpha^\oplus(s);\, s\in \R^+,  f \in C(X) \big\}\\
&=& 
\sup\big\{\int_{Y}sT^+(\frac{h}{s})\, d\nu  -\alpha^\oplus(s)-\int_{X}h\, d\mu;\, s\in \R^+,  h \in C(X) \big\},
\end{eqnarray*}
which means that $\alpha ({\cal T})$ is a forward convex transfer with Kantorovich operators $T_s^+f=sT^+(\frac{h}{s})-\alpha^\oplus(s)$.\\

For 2) use the fact that $(-\beta)$ is convex decreasing to write
\begin{eqnarray*}
{\mathcal T}(\mu, \nu))&=&\sup\{\int_Yg\, d\nu-\beta \big(\int_XT^-g\, d\mu); g\in C(Y)\}\\
&=&
\sup\{\int_Y g\, d\nu+ \sup_{s>0}\{\int_X-sT^-g\, d\mu -(-\beta)^\ominus(s)\}; g\in C(Y)\}\\
&=&
\sup\{\int_Y g\, d\nu- s\int_XT^-g\, d\mu -(-\beta)^\ominus(s); s>0, g\in C(Y)\}.
\end{eqnarray*}
Hence ${\cal T}$ is the supremum of backward linear transfers.\\

\noindent{\it Example 14: General entropic functionals are backward completely convex transfers}

 Consider the following {\it generalized entropy,}
 \begin{equation}
 {\cal T}_\alpha (\mu, \nu)=\int_X \alpha (|\frac{d\nu}{d\mu}|)\,  d\mu, \quad \hbox{if $\nu<<\mu$ and $+\infty$ otherwise,}
 \end{equation}
 where  $\alpha$ is any strictly convex lower semi-continuous superlinear (i.e., $\lim\limits_{t\to +\infty}\frac{\alpha (t)}{t}=+\infty$) real-valued function on $\R^+$. 
 It is then easy to show \cite{GRS}  that 
 \begin{equation}
  {\cal T}_\mu^*(f)=\inf\{\int_X[\alpha^\oplus(f(x)+t)-t]\, d\mu(x); t\in \R\}, 
 \end{equation}
In other words, ${\cal T}_\alpha$ is a backward completely convex transfer, with Kantorovich operators
\[
T_t^-f(x)=\alpha^\ominus(f(x)+t)-t.
\]
  
\noindent{\it Example 15: The logarithmic entropy is a backward $\log$-transfer}

 The relative logarithmic entropy ${\cal H}(\mu, \nu)$ 
  is defined as
\[
\hbox{${\cal H}(\mu, \nu):=\int_X\log (\frac {d\nu}{d\mu})\, d\nu$ if $\nu<<\mu$ and $+\infty$ otherwise.}
\]
It can also be written as 
\[
\hbox{${\cal H}(\mu, \nu):=\int_X h(\frac {d\nu}{d\mu})\, d\mu$ if $\nu<<\mu$ and $+\infty$ otherwise,}
\]
where $h(t)=t\log t -t+1$, which is strictly convex and positive. Since $h^*(t)=e^t-1$, it follows that 
 \[
  {\cal H}_\mu^*(f)=\inf\{\int_X(e^te^{f(x)}-1-t)\, d\mu(x); t\in \R\}=\log\int_Xe^f\, d\mu.
 \]
 In other words, 
${\cal H}(\mu, \nu)=\sup\{\int_X f\, d\nu-\log\int_Xe^f\, d\mu; f\in C(X)\}$, 
 and ${\cal H}$ is therefore  {\it a backward $\beta$-transfer} with $\beta (t)=\log t$, and $E ^-f=e^{f}$ is a convex Kantorovich operator. \\
 ${\cal H}$ is a convex backward transfer since
\begin{eqnarray*}
{\cal H}(\mu, \nu)&=&\sup\{\int_X f\, d\nu-\log\int_Xe^f\, d\mu; f\in C(X)\}\\
&=&
\sup\{\int_X f\, d\nu+ \sup_{s>0}\{\int_X-se^f\, d\mu -\beta^\ominus(s)\}; f\in C(X)\}\\
&=&
\sup\{\int_X f\, d\nu- s\int_Xe^f\, d\mu -\beta^\ominus(s); s>0, f\in C(X)\}.
\end{eqnarray*}
In other words, it is backward completely convex with Kantorovich operators $T^-_sf= se^f+\beta^\ominus(s)$ where $s>0$. \\

\noindent{\it Example 16: A mean-field planning problem} (Orrieri-Porretta-Savar\'e \cite{Sav})

Let $L:\R^d\times \R^d \to \R$ be a Tonelli Lagrangian and $F:\R^d \times L^\infty([0, T];{\cal P}(\R^d)) \to \R$ be a functional that is convex in the second variable and consider the following correlation between two probability measures $\mu$ and $\nu$,    
   \begin{equation}
{\cal T}(\mu, \nu)=:\min
\int_0^T\int_{\R^d} L(x,\vv)\,\rho(t, dx)\, dt+
\int_0^TF(x,\rho(t, dx))\,\d t \,:\quad \vv\in
L^2(\rho (t, dx)\, dt),
\end{equation}
subject to $\rho$ and $v$ satisfying
\begin{eqnarray}
 \label{con} 
\partial_t \rho-\Delta \rho+\nabla\cdot (\rho\,\vv )=0, \quad
\rho(0,\cdot)=\mu\,, \rho(T,\cdot)=\nu.  
\end{eqnarray}
One can show that ${\cal T}$ is  both a completely convex forward and backward transfer. 
Indeed, for each $\ell\in C([0, T], \R^d)$, we consider the Kantorovich operator on $C(\R^d)$, 
\[
T_\ell (u)= u_\ell (T, x)-   \iint_Q F^*(x, \ell(t,x))\,\d x
\]
 where $u_\ell(t, x)$ is a solution of the Hamilton-Jacobi equation 
\begin{eqnarray*}
 -\partial_t u+H(x,Du)&=&\ell\quad\text{in } (0, T)\times \R^d, \\
 u(0, x)&=&u(x).
 \end{eqnarray*}
A standard min-max argument then yields that
\begin{eqnarray*} 
{\cal T}(\mu, \nu)&=&\sup \left\{\int_{\R^d}T_\ell u\, d\nu-\int_{\R^d}u \, d\mu; u\in C(\R^d), \ell \in C([0, T], \R^d)\right\}\\
&=&\sup \left\{\int_{\R^d}u_\ell (T, x) d\nu-
  \int_{\R^d}u_\ell (0, x)\, d\mu (x)-
                     \iint_Q F^*(x, \ell(t,x))\,\d x;\,  
                     -\partial_t u_\ell+H(x,Du_\ell)=\ell \, \, \text{in }Q\right\}.
\end{eqnarray*}
Here $F^*(x, \ell)=\sup \left\{\langle \ell, \rho\rangle -F(x, \rho); m\in L^\infty([0, T];{\cal P}(\R^d))\right\}$ and $Q=(0, 1)\times \R^d$.  \\

Another (stochastic) completely convex -but only backward- transfer can be defined as 
 \begin{equation}
{\cal T}(\mu, \nu)=:\min
\int_0^T\int_{\R^d} L(x,\vv)\,\rho(t, dx)\, dt+
\int_0^TF(x,\rho(t, dx))\,\d t \,:\quad \vv\in
L^2(\rho (t, dx)\, dt),
\end{equation}
subject to $\rho$ and $v$ satisfying
\begin{eqnarray} \label{con.bis} 
\partial_t \rho-\Delta \rho+\nabla\cdot (\rho\,\vv )=0, \quad
\rho(0,\cdot)=\mu\,, \rho(T,\cdot)=\nu.  
\end{eqnarray}

 
\noindent{\it Example 17: The Fisher-Donsker-Varadhan information is a backward completely convex transfer} \cite{DV}

Consider an $\XX$-valued time-continuous Markov process $(\Omega, \FF,
(X_t)_{t\ge0}, (\pp_x)_{x\in \XX})$  with an invariant probability
measure $\mu.$ Assume the transition semigroup, denoted
$(P_t)_{t\ge0},$ to be completely continuous on
$L^2(\mu):=L^2(\XX,\BB,\mu)$. Let $\LL$ be its generator with
domain $\dd_2(\LL)$ on $L^2(\mu)$ and assume the corresponding Dirichlet form
$
\EE(g,g):=\<-\LL g, g\>_{\mu} $ for $g\in \dd_2(\LL)
$
is closable in $L^2(\mu),$ with closure $(\EE, \dd(\EE))$. The Fisher-Donsker-Varadhan information of $\nu$ with
respect to $\mu$ is defined by
\begin{equation}
{\cal I}(\mu|\nu):=\begin{cases}\EE(\sqrt{f}, \sqrt{f}), \ \ &\text{ if }\ \nu=f\mu, \sqrt{f}\in\dd(\EE)\\
+\infty, &\text{ otherwise.}
\end{cases}
\end{equation}
Note that when $(P_t)$ is $\mu$-symmetric, $\nu\mapsto I(\mu|\nu)$ is
exactly the Donsker-Varadhan entropy i.e.\! the rate function
governing the large deviation principle of the empirical measure
$
L_t:=\frac 1t\int_0^t \delta_{X_s} ds
$
for large time $t$.  The corresponding Feynman-Kac semigroup on $L^2(\mu)$
\begin{equation}\label{Feynman-Kac} P_t^u g(x):=\ee^x g(X_t)
\exp\left(\int_0^tu(X_s)\,ds\right).
\end{equation}
It has been proved in \cite{Wu} that ${\cal I}_\mu^*(f)=\log\|P_1^f\|_{L^2(\mu)}$, which yields that $\cal I$ is a backward completely convex transfer.  
\begin{eqnarray*}
{\cal I}_\mu^*(f)=\log\|P_1^f\|_{L^2(\mu)}=\frac{1}{2}\log\|P_1^f\|^2_{L^2(\mu)}=\frac{1}{2}\log \sup\{\int | P_{_1}^fg|^2\, d\mu; \|g\|_{L^2(\mu)}\leq 1\}.
\end{eqnarray*}
In other words, with $\beta (t)=\log t$, we have 
\begin{eqnarray*}
{\cal I}(\mu, \nu)&=&\sup\{\int_Y f\, d\nu-\frac{1}{2}\log \sup\{\int_X | P_1^fg|^2\, d\mu; \|g\|_{L^2(\mu)}\leq 1\}; f\in C(X)\}\\
&=&\sup\{\int_Y f\, d\nu+\sup_{s>0}\sup\limits_{\|g\|_{L^2(\mu)}\leq 1}\frac{1}{2}\{\int_X (-s| P_1^fg|^2-\beta^\ominus(s))\, d\mu  \}; f\in C(X)\}\\
&=&\sup\{\int_Y f\, d\nu-\inf_{s>0}\inf\limits_{\|g\|_{L^2(\mu)}\leq 1}\frac{1}{2}\{\int_X (s| P_1^fg|^2+\beta^\ominus(s))\, d\mu  \}; f\in C(X)\}\\
&=&\sup\{\int_Y f\, d\nu-\int_X T^-_{s, g}f\, d\mu \, ;s \in \R^+, \|g\|_{L^2(\mu)}\leq 1, f\in C(X)\}.
 \end{eqnarray*}
Hence, it is a backward completely convex transfer, with convex Kantorovich operators $(T^-_{s, g})_{s,g}$ defined by $T^-_{s, g}f=\frac{s}{2}| P_{_1}^fg|^2+\frac{1}{2}\beta^\ominus(s)$.\\

  \noindent {\it Example 18: A convex transfer which is not completely convex}
  
Let $\Omega \subset \R^d$ be a Borel measurable subset with $1 < |\Omega| < \infty$, $\lambda := \frac{1}{|\Omega|}$, and define for any two given probability measures $\mu$, $\nu$ on $\Omega$, the correlation,
\begin{equation}
\T_\lambda(\mu,\nu) = 
\begin{cases}
0 & \text{if }\nu \in \mcal{C}_\lambda(\mu)\\
+\infty & \text{otherwise,}
\end{cases}
\end{equation}
where $\mcal{C}_\lambda(\mu) := \{ \nu \in \mcal{P}(\Omega)\,;\, \lambda\lf|\frac{\d\nu}{\d\mu}\rt| \leq 1\, \mu\text{-a.e.}\}$. Note that when $\mu = \lambda \d x|_{\Omega}$ (the uniform measure on $\Omega$), 
\begin{equation}
\T_\lambda(\lambda \d x|_{\Omega},\nu) =
\begin{cases}
0 & \text{if $\lf|\frac{\d\nu}{\d x}\rt| \leq 1$ Lebesgue-a.e.}\\
+\infty & \text{otherwise.}
\end{cases}
\end{equation}
We claim that $\T_\lambda$ is a backward convex transfer but not a \textit{completely} backward convex transfer. Indeed, for the first claim, consider $\alpha_m(t) := (\lambda t)^m\log (\lambda t)$ for $m \geq 1$ and $t \geq 0$, and define
\begin{equation}
\T_m(\mu,\nu) := 
\begin{cases}
\int_{\Omega}\alpha_m\lf(\lf|\frac{\d\nu}{\d\mu}\rt|\rt)\d \mu, & \text{if } \nu << \mu,\\
+\infty & \text{otherwise.}
\end{cases}
\end{equation}
By Example 14, $\T_m$ is a backward completely convex transfer, and 
\begin{equation}
(\T_{m,\mu})^*(f) = \inf\{\int_{\Omega}[\alpha_m^{\oplus}(f(x) + t) - t]\d\mu(x)\,;\, t \in \R\}.
\end{equation}
The function $\alpha_m^{\oplus}$ can be explicitly computed as 
\begin{equation}
\alpha^{\oplus}_m(t)= 
\begin{cases} e^{-1 + \frac{1}{m-1} W\lf(\beta_m t \rt) }\lf[ \beta_m t + \frac{1}{m} e^{W\lf(\beta_m t \rt)}\rt] & \text{if } t \geq -\frac{\lambda}{m-1} e^{-1},\\
0 & \text{if }  t < -\frac{\lambda}{m-1} e^{-1}.
\end{cases}
\end{equation}
where $\beta_m := \frac{m-1}{\lambda m}e^{\frac{m-1}{m}}$, and $W$ is the \textit{Lambert-W} function. It is easy to see that $\T_\lambda(\mu,\nu) = \sup_{m}\T_m(\mu,\nu)$; hence it is a backward convex transfer (as a supremum of backward convex transfers). 

However, $\T_\lambda$ is not a \textit{completely} backward convex transfer, since 
$$(\T_{\lambda, \mu})^*(f) = (\sup_{m} \T_{m,\mu})^*(f) \leq \inf_{m} \T_{m,\mu}^*(f) = \int \frac{f}{\lambda}\d \mu,$$
with the inequality being in general strict.

Note that this also implies that the Wasserstein projection on the set ${\cal C}_\mu$, that is 
\begin{equation}
W_2^2(P_1[\nu], \nu ) =\inf\{W_2^2(\sigma, \nu); \lf|\frac{\d\sigma}{\d x}\rt| \leq 1\}= \inf\{ \T_\lambda(\lambda\d x|_{\Omega}, \sigma) + W_2^2(\sigma, \nu)\,;\, \sigma \in \mcal{P}(\Omega)\}
\end{equation}
is in fact an inf-convolution of a convex-\textit{but not completely convex}- transfer $\T_\lambda$ with the linear transfer $W_2^2$, and no duality formula can then be extracted.

\section{Operations on convex and entropic transfers}

Denote by ${\cal CT}_-(X\times Y)$ (resp., ${\cal CT}_+(X\times Y)$) the class of backward (resp., forward) completely convex transfers. They are clearly convex cones in the space of  convex weak$^*$-lower semi-continuous functions on ${\cal P}(X)\times {\cal P}(Y)$. They also satisfy the following permanence properties. The most important being that the inf-convolution with linear transfers  generate many new  examples of convex and entropic transfers..

\begin{prop}\label{con.con} Let ${\mathcal F}$ be a backward  completely convex transfer with  Kantorovich operators $(F)_i ^-$, Then,
 \begin{enumerate}
 
 \item If $a\in \R^+\setminus \{0\}$, then  $a{\mathcal F}\in {\cal CT}_-(X\times Y)$  with Kantorovich operators given by 
$F _{a, i}^-(f)=aF_i ^-(\frac{f}{a}).$
 
\item If %
${\cal T}$ is a backward linear transport on $Y\times Z$ with Kantorovich operator $T^-$,  then ${\mathcal F}\star{\mathcal T}$ is a backward completely convex transfer with Kantorovich operators given by $F _{i}^-\circ T ^-$. 
\end{enumerate}
 \end{prop} 
\noindent{\bf Proof:} Immediate. For 2) we calculate the Legendre dual of $({\mathcal F}\star T)_\mu$ at $g\in C(Z)$ and obtain,
  \begin{eqnarray*}
( {\mathcal F}\star T)_\mu^*(g)&=&\sup\limits_{\nu \in {\mathcal P}(Z)}\sup\limits_{\sigma \in  {\mathcal P}(Y)} 
 \left\{\int_{Z} g\, d\nu -  {\mathcal F}(\mu, \sigma) -{\mathcal T}(\sigma, \nu) \right\}\\
&=&  \sup\limits_{\sigma \in  {\mathcal P}(Y)} 
 \left\{{\mathcal T}_\sigma^* (g)-{\mathcal F}(\mu, \sigma) \right\}\\
&=& \sup\limits_{\sigma \in  {\mathcal P}(Y)} 
 \left\{\int_{Y} T^-g\, d\sigma-{\mathcal F}(\mu, \sigma) \right\}\\
&=&({\mathcal F})_\mu^* (T^-(g))\\
 &=& \inf\limits_{i\in I}\int_{X}F_i^-\circ T^-g(x))\, d\mu(x). 
 \end{eqnarray*} 
 The same properties hold for entropic transfers.
That  we will denote by ${\cal E}$ as opposed to ${\cal T}$ to distinguish them from the linear transfers. We shall use $E^+$ and $E^-$ for their Kantorovich operators.

 \begin{prop}\rm  Let $\beta: \R\to \R$ be a concave increasing function and let ${\mathcal E}$ be a $\beta$-backward  transfer with Kantorovich operator $E ^-$. Then,
 \begin{enumerate}

\item If $\lambda \in \R^+\setminus \{0\}$, then  $\lambda{\mathcal E}$  is a $(\lambda \beta)$-backward transfer   with Kantorovich operator $E _\lambda^-(f)=  E ^-(\frac{f}{\lambda})$.

\item  $\tilde {\mathcal E}$ is a $((-\beta)^\ominus)^\oplus$-forward convex transfer with Kantorovich operator ${\tilde E}^+h=-E^-(-h)$.

\item If ${\mathcal T}$  is a backward linear transfer  on  $Y\times Z$  with Kantorovich operator $T^-$, then ${\mathcal E}\star{\mathcal T}$ is a  a backward  $\beta$-transfer on  $X\times Z$ with Kantorovich operator equal to $E^-\circ T^-$. In other words,  
 \begin{equation}
{\mathcal E}\star{\mathcal T}\, (\mu, \nu)= 
\sup\big\{\int_{Z}g(y)\, d\nu(y)-\beta (\int_{X}E^-\circ T^-g(x))\, d\mu(x));\, g\in C(Z) \big\}. 
\end{equation}

\end{enumerate}
\end{prop}
\noindent {\bf Proof:} 1) is trivial. For 2) note that since $\beta$ is concave and increasing, then 
\begin{eqnarray*}
\tilde {\mathcal T}(\nu, \mu))&=&{\mathcal T}(\mu, \nu))\\
&=&\sup\{\int_Yg\, d\nu-\beta \big(\int_XT^-g\, d\mu); g\in C(Y)\}\\
&=&
\sup\{\int_Y g\, d\nu+ \sup_{s>0}\{\int_X-sT^-g\, d\mu -(-\beta)^\ominus(s)\}; g\in C(X)\}\\
&=&
\sup\{\int_Y g\, d\nu- s\int_XT^-g\, d\mu -(-\beta)^\ominus(s); s>0, g\in C(X)\}\\
&=&
\sup\{s\int_X-T^-(-h)\, d\mu -(-\beta)^\ominus(s)-\int_Y h\, d\nu; s>0, g\in C(X)\\
&=&
\sup\{((-\beta)^\ominus)^\oplus(\int_X-T^-(-h)\, d\mu)-\int_Y h\, d\nu; s>0, h\in C(X)\}.
\end{eqnarray*}
In other words, $\tilde {\mathcal T}$ is a $(\beta^\ominus)^\oplus$-forward convex transfer. \\
For 3) we  calculate the Legendre dual of $({\mathcal E}\star T)_\mu$ at $g\in C(Z)$ and obtain,
 \begin{eqnarray*}
( {\mathcal E}\star T)_\mu^*(g)&=&\sup\limits_{\nu \in {\mathcal P}(Z)}\sup\limits_{\sigma \in  {\mathcal P}(Y)} 
 \left\{\int_{Z} g\, d\nu -  {\mathcal E}(\mu, \sigma) -{\mathcal T}(\sigma, \nu) \right\}\\
&=&  
\sup\limits_{\sigma \in  {\mathcal P}(Y)} 
 \left\{{\mathcal T}_\sigma^* (g)-{\mathcal E}(\mu, \sigma) \right\}\\
&=& \sup\limits_{\sigma \in  {\mathcal P}(Y)} 
 \left\{\int_{Y} T^-g\, d\sigma-{\mathcal E}(\mu, \sigma) \right\}\\
&=&({\mathcal E})_\mu^* (T^-(g))\\
 &=& \beta \Big(\int_{X}E^-\circ T^-g(x))\, d\mu(x)). 
 \end{eqnarray*} 
A similar statement holds for $\alpha$-forward transfers where $\alpha$ is now a convex increasing function on $\R^+$. But we then have to reverse the orders. For example, if ${\mathcal T}$ (resp., ${\mathcal E}$) is a forward linear transfer  on  $Z\times X$ (resp., a forward  $\alpha$-transfer on $X\times Y$) with Kantorovich operator $T^+$ (resp., $E^+$), then ${\mathcal T}\star{\mathcal E}$ is a forward  $\alpha$-transfer on  $Z\times Y$ with Kantorovich operator equal to $E^+\circ T^+$. In other words, 
 \begin{equation}\label{E1}
{\mathcal T}\star{\mathcal E}\, (\mu, \nu)= 
\sup\big\{\alpha\big(\int_{Y}E^+\circ T^+f(y))\, d\nu(y)\big)-\int_{X}f(x)\, d\mu(x);\,  f\in C(X) \big\}.
\end{equation}

\section{Subdifferentials of linear and convex transfers}

If  $ \mathcal{T}$ is a linear transfer, then both $ \mathcal{T}_\mu$ and $ \mathcal{T}_\nu$ are convex weak$^*$ lower semi-continuous and one can therefore consider their (weak$^*$) subdifferential $\partial \T_{\mu}$ (resp., $\partial \T_{\nu}$) in the sense of convex analysis. In other words, 
\[
\hbox{$g \in \partial \T_{\mu}(\nu)$ if and only if  $\T(\mu,\nu') \geq \T(\mu,\nu) + \int_Y g\d (\nu' - \nu)$ \,\, for any $\nu'\in {\mathcal P}(Y)$.}
\]
In other words, $g \in \partial \T_{\mu}(\nu)$ if and only if $\T_{\mu}(\nu) + \T^*_\mu (g) = \langle g, \nu\rangle.$ Since $\T_{\mu}(\nu) = \T(\mu,\nu)$ and $\T^*_\mu (g) = \int T^{-}g \d\mu$, we then obtain the following characterization of the subdifferentials. 

 \begin{prop}Let $\T$ be a backward (resp., forward) linear transfer. Then the subdifferential of $\T_\mu :\mathcal{P}(Y) \to \R\cup\{+\infty\}$ at $\nu \in \mathcal{P}(Y)$ (resp., $\T_\nu:\mathcal{P}(X) \to \R\cup\{+\infty\}$ at $\mu \in \mathcal{P}(X)$) is given by
\begin{equation}
\partial \T_{\mu}(\nu) = \left\{g \in C(Y)\,:\, \int_{Y}g(y)\d\nu(y) - \int_{X}T^- g(x)\d\mu(x) = \T(\mu,\nu)\right\}
\end{equation}
respectively,
\begin{equation}
\partial \T_{\nu}(\mu) = \lf\{f \in C(X)\,:\, \int_{Y}T^{+}f(y)\d\nu(y) - \int_{X} f(x)\d\mu(x) = \T(\mu,\nu)\rt\}
\end{equation}
In other words, the subdifferential of $\T_\mu$ at $\nu$ (resp., $\T_\nu$ at $\mu$) is exactly the set of maximisers for the dual formulation of $\T(\mu,\nu)$. 
\end{prop}
It is easy to see that the same expressions hold - with the necessary modifications - for backward completely convex transfers (resp., forward completely convex transfers), as well as $\beta$-backward transfers (resp., $\alpha$-forward transfers).

In the following, we observe some elementary consequences for elements in the subdifferential.

 \begin{prop} Suppose $\T$ is a linear backward transfer such that the Dirac masses are contained in $D_1(\T)$. Fix $\mu \in \mcal{P}(X)$ and $\nu \in \mcal{P}(Y)$, and suppose the infimum in $\T(\mu,\nu)$ is achieved by $\bar{\pi}$ with disintegration w.r.t. $\mu$ denoted by $\bar{\pi}_x$, that is 
\[
\T(\mu,\nu) = \int_{X}\T(x,\pi_x)\d\mu(x). 
\] 
Then, for each $\bar{f} \in \partial \T_\mu(\nu)$, we have
\eqs{
T^{-}\bar{f}(x) = \int_{Y}\bar{f}(y)\d\bar{\pi}_x(y) - \T(x,\bar{\pi}_x), \quad \text{for $\mu$-a.e. $x \in X$.
}}
 \end{prop}
\rm {\bf Proof:} Indeed, if $\bar{f} \in \partial \T_\mu(\nu)$, then by definition
\eqs{\label{equality}
\int_{Y}\bar{f}(y)\d\nu(y) - \int_{X}T^- \bar{f}(x)\d\mu(x) = \T(\mu,\nu) = \int_{X}\T(x,\bar{\pi}_x)\d\mu(x), 
}
that is
$
\int_{X}\lf[T^- \bar{f}(x) - \int_{Y}\bar{f}(y)\d\bar{\pi}_x(y) + \T(x,\bar{\pi}_x)   \rt]\d\mu = 0.
$
Since $T^- \bar{f}(x) = \sup_{\sigma}\lf\{\int \bar{f}\d\sigma - \T(x,\sigma)\rt\}$, the quantity in the brackets  is non-negative and we get our claim.
 \begin{prop}
Suppose $\T$ is a linear backward transfer such that the Dirac masses are contained in $D_1(\T)$. Fix $\mu$, and suppose $\nu \mapsto \T(\mu,\nu)$ is strictly convex. If $\partial \T_\mu(\nu)$ is non-empty for some $\nu \in \mcal{P}(Y)$, then the infimum in $\T(\mu,\nu)$ is attained, i.e., 
there exists $\bar{\pi}$ such that $\T(\mu,\nu)=\int_{X}\T(x,\pi_x)\d\mu(x)$.
 \end{prop}
{\bf Proof:} If  $\bar{f}\in \partial \T_\mu(\nu)$, then 
$
\int \bar{f}\d\nu - \int T^{-}\bar{f}\d\mu = \T(\mu,\nu). $
From the expression $T^- \bar{f}(x) = \sup_{\sigma}\lf\{\int \bar{f}\d\sigma - \T(\delta_x,\sigma)\rt\}$, we know the supremum will be achieved by some $\sigma_x$. Defining $\tilde{\pi}$ by $\d\tilde{\pi}(x,y) = \d\mu(x)\d\sigma_x(y)$, and the right marginal of $\tilde{\pi}$ by $\tilde{\nu}$, we integrate against $\mu$ to achieve
\eqs{
\int T^- \bar{f}\d\mu = \int \bar{f}\d\tilde{\nu} - \int \T(\delta_x, \sigma_x)\d\mu.
}
This shows that $\T(\mu,\tilde{\nu}) = \inf_{\pi \in \Gamma(\mu,\tilde{\nu})}\int \T(\delta_x,\pi_x)\d\mu = \int \T(\delta_x, \sigma_x)\d\mu$, 
and consequently, $\bar{f} \in \partial \T_\mu(\tilde{\nu})$. But by strict convexity, this can only be true if $\tilde{\nu} = \nu$.

The following can be seen as Euler-Lagrange equations for variational problems on spaces of measures, and follows closely \cite{FGJ}.

\begin{prop} Let $\T_\alpha(\mu,\nu) := \int_{X} \alpha\lf(\frac{\d\nu}{\d\mu}\rt)\d\mu$ be the generalised entropy transfer considered in Example 14, and let $\T$ be any linear backward transfer. For a fixed $\mu$, consider  the functional  
$
I_\mu(\nu) := \T_\alpha(\mu,\nu) - \T(\mu,\nu),
$
and assume $\bar{\nu}$ realises $\inf_{\nu \in \mcal{P}(X)}I_\mu(\nu)$. Then, there exists $\bar{f} \in \partial \T_\mu(\bar{\nu})$ such that the following Euler-Lagrange equation holds for $\bar{\nu}-$a.e. $x \in X$,
\eqs{
\alpha'\lf(\frac{\d\bar{\nu}}{\d\mu}\rt) = \bar{f} + C, 
}
where $C$ is a constant. 

If $\T_\alpha$ is replaced with the logarithmic entropic transfer $\mcal{H}(\mu,\nu) = \int \log(\frac{\d\nu}{\d\mu})\d\nu$, then 
\eqs{
\log \lf(\frac{\d\bar{\nu}}{\d\mu}\rt) = \bar{f} + C.  }

\end{prop}
{\bf Proof:} Recall that $\T_\alpha(\mu,\nu) := \int_{X}\alpha(|\frac{\d\nu}{\d\mu}|) \d\mu$ if $\nu << \mu$ (and $+\infty$ otherwise) is a backward completely convex transfer, with
\eqs{
\T_{\mu}^*(f) = \inf\lf\{\int_{X}[\alpha^{\oplus}(f(x) + t) - t]\d\mu(x)\,;\, t \in \R\rt\},
}
where 
$
T_t^- f(x) := \alpha^{\ominus}(f(x) + t) - t.
$
are the corresponding Kantorovich transfers. Here $\alpha \in C^1$, is strictly convex and superlinear. It follows that  
\eqs{
\alpha'(|\frac{\d\nu}{\d\mu}|) \in \partial \T_{\mu}(\nu).
}
We can see this either directly from the subdifferential definition, or from observing 
$$\alpha^{\oplus}(\alpha'(|\frac{\d\nu}{\d\mu}|)) = \frac{\d\nu}{\d\mu}\alpha'(|\frac{\d\nu}{\d\mu}|) - \alpha(|\frac{\d\nu}{\d\mu}|).$$
In particular,
\as{
\T_{\mu}^*\lf(\alpha'(|\frac{\d\nu}{\d\mu}|)\rt) = \int_{X}\alpha^{\oplus}\lf(\alpha'(|\frac{\d\nu}{\d\mu}|)\rt) \d\mu .
}
The rest is an easy adaptation of Theorem 2.2 in \cite{FGJ}.

\section{Transfer inequalities}
Let ${\mathcal T}$ be a transfer, and let ${\mathcal E}_1$, ${\mathcal E}_2$ be entropic transfers on $X\times X$. Standard Transport-Entropy or Transport-Information inequalities are usually of the form 
\begin{equation}\label{one}	
\hbox{${\mathcal T}(\sigma, \mu) \leq \lambda_1 {\mathcal E}_1(\mu, \sigma)$ \quad for all $\sigma \in {\cal P}(X)$,} 
\end{equation}
\begin{equation}\label{two}	
\hbox{${\mathcal T}(\mu, \sigma) \leq \lambda_2 {\mathcal E}_2(\mu, \sigma)$ \quad for all $\sigma \in {\cal P}(X)$,} 
\end{equation}
\begin{equation}\label{three}
\hbox{${\mathcal T}(\sigma_1, \sigma_2) \leq \lambda_1 {\mathcal E}_1(\sigma_1, \mu)+\lambda_2 {\mathcal E}_2( \sigma_2, \mu)$ \quad for all $\sigma_1, \sigma_2 \in {\cal P}(X)$,} 
\end{equation}
where $\mu$ is a fixed measure, and $\lambda_1$, $\lambda_2$ are two positive reals.
In our terminology, Problem \ref{one} (resp., \ref{two}), (resp., \ref{three}) amount to find $\mu$,  $\lambda_1$, and $\lambda_2$ such that 
\begin{equation}
(\lambda_1{\mathcal E}_1)\star (-{\mathcal T}) \, (\mu, \mu)\geq 0,
\end{equation} 
\begin{equation}(\lambda_2{\mathcal E}_2)\star (-\tilde {\mathcal T})\,  (\mu, \mu)\geq 0,
\end{equation} 
 \begin{equation}(\lambda_1\tilde {\mathcal E}_1)\star (-{\mathcal T})\star (\lambda_2{\mathcal E}_2) \, (\mu, \mu)\geq 0,
\end{equation} 
where $\tilde {\mathcal T} (\mu, \nu)={\mathcal T} (\nu, \mu)$. Note for example that 
\[
\tilde {\mathcal E}_1\star  (-{\mathcal T})\star  {\mathcal E}_2\, (\mu, \nu)= \inf\{\tilde{\mathcal E}_1(\mu, \sigma_1)-{\mathcal T}_2(\sigma_1, \sigma_2)+{\mathcal E}_2(\sigma_2, \nu);\, \sigma_1, \sigma_2 \in {\mathcal P}(Z)\}.
\]
We shall therefore write duality formulas for the transfers ${\mathcal E}_1\star (-{\mathcal T})$, ${\mathcal E}_2\star (-\tilde {\mathcal T})$ and $\tilde {\mathcal E}_1\star (-{\mathcal T})\star {\mathcal E}_2$ between any two measures $\mu$ and $\nu$, where ${\mathcal T}$ is any convex transfer, while ${\mathcal E}_1$, ${\mathcal E}_2$ are entropic transfers.

\subsection{Backward convex to backward completely convex inequalities}

 We would like to prove inequalities such as 
\begin{equation}\label{one.main}	
\hbox{${\mathcal F}_2(\sigma, \mu) \leq {\mathcal F}_1(\mu, \sigma)$ \quad for all $\sigma \in {\cal P}(X)$,} 
\end{equation}
where both ${\mathcal F}_1$ and ${\mathcal F}_2$ are backward convex transfers. We then apply it to Transport-Entropy inequalities of the form 
\begin{equation}\label{one.revised}	
\hbox{${\mathcal F}(\sigma, \mu) \leq \lambda {\mathcal E}\star {\cal T}(\mu, \sigma)$ \quad for all $\sigma \in {\cal P}(X)$,} 
\end{equation}
where ${\mathcal F}$ is a backward convex transfer, while ${\mathcal E}$ is a $\beta$-entropic transfer and ${\cal T}$ is a backward linear transfer. 

\begin{prop} \label{back-back} 
Let ${\mathcal F}_1$ be a backward completely convex transfer with Kantorovich operator $(F^-_{1,i})_{i\in I}$ on $X_1\times X_2$, and ${\mathcal F}_2$ is a backward convex transfer on $X_2\times X_3$ with Kantorovich operator $(F^-_{2, j})_{j\in J}$. 
\begin{enumerate}

\item The following duality formula hold: 
\begin{equation}
{\mathcal F}_1\star -{\mathcal F}_2\, (\mu, \nu)=\inf\limits_{f\in C(X_3)}\inf\limits_{j\in J}\sup\limits_{i\in I}\left\{-\int_{X_1} F_{1,i}^-\circ F _{2,j}^-f\, d\mu-\int_{X_3} f\, d\nu\right\}.
\end{equation}
\item If ${\mathcal F}_1$ is a $\beta$-backward transfer on $X_1\times X_2$ with Kantorovich operator $E_1^-$, 
 then
\begin{equation}
{\mathcal F}_1\star -{\mathcal F}_2\, (\mu, \nu)=\inf\limits_{f\in C(X_3)}\inf\limits_{j\in J}\left\{-\beta (\int_{X_1} E_1^-\circ F_{2,j}^-f\, d\mu)-\int_{X_3} f\, d\nu  
\right\}.
\end{equation}
 \end{enumerate}
\end{prop}
\noindent{\bf Proof:} Write
\begin{eqnarray*}
{\mathcal F}_1\star -{\mathcal F}_2\, (\mu, \nu)&=& \inf\{{\mathcal F}_1(\mu, \sigma)-{\mathcal F}_2(\sigma, \nu);\, \sigma \in {\mathcal P}({X_2}) \}\\
 &=&
\inf\limits_{\sigma \in {\mathcal P}({X_2})}\left\{{\mathcal F}_1(\mu, \sigma)- \sup\limits_{f\in C(X_3)}\sup\limits_{j\in J}\left\{\int_{X_3} f\, d\nu-\int_{X_2} F _{2,j}^-f\, d\sigma\right\}\right\}\\
&=&
\inf\limits_{\sigma \in {\mathcal P}({X_2})}\inf\limits_{f\in C(X_3)}\inf\limits_{j\in J}\left\{{\mathcal F}_1(\mu, \sigma)- \int_{X_3} f\, d\nu+\int_{X_2} F _{2,j}^-f\, d\sigma \right\}\\
&=&\inf\limits_{f\in C(X_3)}\inf\limits_{j\in J}\left\{-\sup\limits_{\sigma \in {\mathcal P}({X_2})}\{-\int_{X_2} F _{2,j}^-f\, d\sigma    -{\mathcal F}_1(\mu, \sigma)\}- \int_{X_3} f\, d\nu  \right\}\\
&=&\inf\limits_{f\in C(X_3)}\inf\limits_{j\in J}\left\{-({\mathcal F}_{1})_\mu^*(-F _{2,j}^-f)-\int_{X_3} f\, d\nu  
\right\}\\
&=&\inf\limits_{f\in C(X_3)}\inf\limits_{j\in J}\left\{-\inf\limits_{i\in I}\int_{X_1} F_{1,i}^-\circ -F _{2,j}^-f\, d\mu-\int_{X_3} f\, d\nu\right\}\\
&=&\inf\limits_{f\in C(X_3)}\inf\limits_{j\in J}\sup\limits_{i\in I}\left\{-\int_{X_1} F_{1,i}^-\circ -F _{2,j}^-f\, d\mu-\int_{X_3} f\, d\nu\right\}.
\end{eqnarray*}
2) If ${\mathcal F}_1$ is a $\beta$-backward transfer on $X_1\times X_2$ with Kantorovich operator $E_1^-$, 
 then use in the above calculation that $({\mathcal F}_{1})_\mu^*(-F _{2,j}^-f)=\beta (\int_{X_1} E_1^-\circ -F_{2,j}^-f\, d\mu)$.

\begin{cor} Let  ${\mathcal F}$ be a backward convex transfer on $Y_2\times X_2$ with Kantorovich operators $(F^-_i)_{i\in I}$ and let ${\mathcal E}$ be a backward $\beta$-transfer on $X_1\times Y_1$ with Kantorovich operator $E^-$. Let ${\cal T}$ be a backward linear transfer  on $Y_1\times Y_2$ with Kantorovich operator $T^-$ and $\lambda>0$. Then, for any fixed pair of probability measures $\mu\in {\cal P}( X_1)$ and $\nu\in {\cal P}(X_2 )$, the following are equivalent:
\begin{enumerate}
\item For all $\sigma \in {\cal P}(Y_2)$, we have ${\mathcal F}(\sigma, \nu) \leq \lambda \, {\mathcal E}\star {\mathcal T}\, (\mu, \sigma)$.
 
\item  For all $g\in C(X_2)$ and $i\in I$, we have $\beta \big( \int_{X_1} E^-\circ T^-\circ \frac{-1}{\lambda}F_i^-(\lambda g)\, d\mu\big)+\int_{X_2} g\, d\nu \leq 0$. 
 \end{enumerate}
\end{cor}
In particular, if we apply the above in the case where $\cal E$ is the logarithmic entropy, that is 
  \begin{equation}
\hbox{
${\cal H}(\mu, \nu)=\int_X\log (\frac {d\nu}{d\mu})\, d\nu$ if $\nu<<\mu$ and $+\infty$ otherwise,}
\end{equation}
 which is a backward $\beta$-transfer with $\beta (t)=\log t$ and $E ^-f=e^f$ as a backward Kantorovich operator.
 
  \begin{cor} Let  ${\mathcal F}$ be a backward convex transfer on $X_2\times Y_2$ with Kantorovich operators $(F^-_i)_{i\in I}$ and let ${\mathcal E}$ be a backward $\beta$-transfer on $X_1\times Y_1$. with Kantorovich operator $E^-$. Let ${\cal T}$ be a backward linear transfer  on $Y_1\times Y_2$ with Kantorovich operator $T^-$ and $\lambda>0$. Then, for any fixed pair of probability measures $\mu\in {\cal P}( X_1)$ and $\nu\in {\cal P}(X_2 )$, the following are equivalent:
 
 \begin{enumerate}
 \item  For all $\sigma \in {\cal P}(Y)$, we have ${\mathcal F}(\sigma, \nu) \leq \lambda\,  {\mathcal H}\star {\mathcal T}\, (\mu, \sigma)$.

\item For all $g\in C(X_2)$, we have  $\sup\limits_{i\in I} \int_{X_1} e^{T^-\circ \frac{-1}{\lambda}F_i^-(\lambda g)}\, d\mu \leq e^{-\int_{X_2} g\, d\nu}$.
\end{enumerate}
In particular, if  ${\cal T}$ is the identity transfer and ${\cal F}$ is a backward linear transfer, then 
  \begin{equation}
\hbox{${\mathcal F}(\sigma, \nu) \leq \lambda\,   {\mathcal H}\, (\sigma, \mu)$ for all $\sigma \in {\cal P}(Y)$  \quad $\Leftrightarrow$ \quad $ \int_{X_1} e^{-F^-(\lambda g)}\, d\mu \leq e^{-\frac{1}{\lambda}}e^{-\int_{X_2} g\, d\nu}$ \quad  for all $g\in C(X_2)$. }
  \end{equation}
  \end{cor}

\subsection{Forward convex to backward completely convex transfer inequalities}

 We are now interested in inequalities such as 
\begin{equation}\label{one.main}	
\hbox{${\mathcal F}_2(\nu, \sigma) \leq {\mathcal F}_1(\mu, \sigma)$ \quad for all $\sigma \in {\cal P}(X)$,} 
\end{equation}
where both ${\mathcal F}_1$ and ${\mathcal F}_2$ are backward convex transfers, and in particular,  Transport-Entropy inequalities of the form 
\begin{equation}\label{one.revised}	
\hbox{${\mathcal F}(\nu, \sigma) \leq \lambda {\mathcal E}\star {\cal T}(\mu, \sigma)$ \quad for all $\sigma \in {\cal P}(X)$,} 
\end{equation}
where  
${\mathcal E}$ is a $\beta$-entropic transfer and ${\cal T}$ is a backward linear transfer. But we can write (\ref{one.main}) as 
\begin{equation}\label{one.main}	
\hbox{$\tilde {\mathcal F}_2(\sigma, \nu) \leq {\mathcal F}_1(\mu, \sigma)$ \quad for all $\sigma \in {\cal P}(X)$,} 
\end{equation}
where now $\tilde {\mathcal F}_2(\sigma, \nu)={\mathcal F}_2(\nu, \sigma)$ is a forward convex transfer. So, we need to establish the following type of duality.

\begin{prop} \label{forward-back}  
\label{B-F} Let ${\mathcal F}_1$ be a backward completely convex transfer with Kantorovich operator $(F^-_{1,i})_{i\in I}$ on $X_1\times X_2$, and let ${\mathcal F}_2$ be a forward convex transfer on $X_2\times X_3$ with Kantorovich operator $(F^+_{2, j})_{j\in J}$. 
\begin{enumerate}

\item The following duality formula then holds: 
\begin{equation}
{\mathcal F}_1\star -{\mathcal F}_2\, (\mu, \nu)=\inf\limits_{g\in C(X_2)}\inf\limits_{j\in J}\sup\limits_{i\in I}\left\{-\int_{X_1} F_{1,i}^-(-g)\, d\nu-\int_{X_3} F _{2,j}^+(g)\, d\nu\right\}.
\end{equation}
\item If ${\mathcal F}_1$ is a $\beta$-backward transfer on $X_1\times X_2$ with Kantorovich operator $E_1^-$, 
 then
\begin{equation}
{\mathcal F}_1\star -{\mathcal F}_2\, (\mu, \nu)=\inf\limits_{g\in C(X_2)}\inf\limits_{j\in J}\left\{-\beta (\int_{X_1} E_1^-(-g)\, d\mu)-\int_{X_3} F _{j}^+(g)\, d\nu\right\}.
\end{equation}
\item If ${\mathcal F}_1$ is a backward $\beta$-transfer with Kantorovich operator $E ^-_1$, and ${\mathcal F}_2$ is a forward $\alpha$-transfer with Kantorovich operator $E ^+_2$, then
\begin{equation}
{\mathcal F}_1\star -{\mathcal F}_2\, (\mu, \nu)=\inf\limits_{g\in C(X_2)}\left\{ -\beta (\int_{X_1} E_1^-(-g)\, d\mu)-\alpha (\int_{X_3} E_2^+g\, d\nu)\right\}.
\end{equation}
\item In particular, if ${\mathcal E}$ is a backward $\beta$-transfer with Kantorovich operator $E ^-$, and ${\mathcal T}$ is a forward linear transfer with Kantorovich operator $T ^+$, then
\begin{equation}
{\mathcal E}\star -{\mathcal T}\, (\mu, \nu)=\inf\limits_{g\in C(X_2)}\left\{ -\beta (\int_{X_1} E^-(-g)\, d\mu)-\int_{X_3} T^+g\, d\nu\right\}.
\end{equation}

\end{enumerate}
\end{prop}
\noindent{\bf Proof:} 1) Assume ${\mathcal F}_1$ is a backward completely convex transfer with Kantorovich operator $F ^-_{1,i}$, and ${\mathcal F}_2$ is a forward convex transfer with Kantorovich operator $F^+_{2, j}$, then
\begin{eqnarray*}
{\mathcal F}_1\star -{\mathcal F}_2\, (\mu, \nu)&=& \inf\{{\mathcal F}_1(\mu, \sigma)-{\mathcal F}_2(\sigma, \nu);\, \sigma \in {\mathcal P}({X_2}) \}\\
&=&
\inf\limits_{\sigma \in {\mathcal P}({X_2})}\left\{{\mathcal F}_1(\mu, \sigma)- \sup\limits_{g\in C(X_2)}\left\{\sup\limits_{j\in J}(\int_{X_3} F _{2,j}^+g\, d\nu) -\int_{X_2} g\, d\sigma\right\}\right\}\\
&=&
\inf\limits_{\sigma \in {\mathcal P}({X_2})}\inf\limits_{g\in C(X_2)}\inf\limits_{j\in J}\left\{{\mathcal F}_1(\mu, \sigma)- \int_{X_3} F _{2,j}^+g\, d\nu +\int_{X_2} g\, d\sigma \right\}\\
&=&\inf\limits_{g\in C(X_2)}\inf\limits_{j\in J}\left\{-\sup\limits_{\sigma \in {\mathcal P}({X_2})}\{-\int_{X_2} g\, d\sigma   -{\mathcal F}_1(\mu, \sigma)\}- \int_{X_3} F _{2,j}^+g\, d\nu)  \right\}\\
&=&\inf\limits_{g\in C(X_2)}\inf\limits_{j\in J}\left\{-({\mathcal F}_{1})_\mu^*(-g)-\int_{X_3} F _{2,j}^+(g)\, d\nu\right\}\\
&=&\inf\limits_{g\in C(X_2)}\inf\limits_{j\in J}\left\{-(\inf\limits_{i\in I}\int_{X_1} F_{1,i}^-(-g)\, d\nu-\int_{X_3} F _{2,j}^+(g)\, d\nu\right\}\\
&=&\inf\limits_{g\in C(X_2)}\inf\limits_{j\in J}\sup\limits_{i\in I}\left\{-\int_{X_1} F_{1,i}^-(-g)\, d\nu-\int_{X_3} F _{2,j}^+(g)\, d\nu\right\}
\end{eqnarray*}
2) If  ${\mathcal F}_1$ is a $\beta$-backward entropic transfer with Kantorovich operator $E ^-$, it suffices to note in the above proof that $({\mathcal F_1})_\mu^*(g)=\beta (\int_X E_1^-(-g)\, d\mu).$\\
 3) If now  ${\mathcal F}_2$ is a forward $\alpha$-transfer with Kantorovich operator $E ^+_2$, then it suffices to note in the above proof that $({\mathcal F}_2)_\nu^*(g)=\alpha (\int_X E_2^+g\, d\nu).$\\
 4) corresponds to when $\alpha (t)=t$.
 
\begin{cor} Let  ${\mathcal F}$ be a backward convex transfer on $X_2\times Y_2$ with Kantorovich operators $(F^-_i)_{i\in I}$ and let ${\mathcal E}$ be a backward $\beta$-transfer on $X_1\times Y_1$ with Kantorovich operator $E^-$.  Let ${\cal T}$ be a backward linear transfer  on $Y_1\times Y_2$ with Kantorovich operator $T^-$ and $\lambda>0$. Then, for any fixed pair of probability measures $\mu\in {\cal P}( X_1)$ and $\nu\in {\cal P}(X_2 )$, the following are equivalent:
\begin{enumerate}
\item For all $\sigma \in {\cal P}(Y_2)$, we have ${\mathcal F}(\nu, \sigma) \leq \lambda \, {\mathcal E}\star {\mathcal T}\, (\mu, \sigma)$. 
 \item  For all $g\in C(X_2)$, we have $\beta  \big( \int_{X_1} E^-\circ T^-g)\, d\mu\big)\leq \inf\limits_{i\in I}\frac{1}{\lambda}\int_{X_2} F_i^-(\lambda g) d\nu$.
 \end{enumerate}
 In particular, if ${\mathcal E}_2$ is a backward $\beta_2$-transfer on $X_2\times Y_2$ with Kantorovich operator $E^-_2$, and ${\mathcal E}_1$ is a backward $\beta_1$-transfer on $X_1\times Y_1$ with Kantorovich operator $E^-_1$, then the following are equivalent:
 \begin{enumerate}
\item  For all $\sigma \in {\cal P}(Y_2)$, we have ${\mathcal E}_2(\nu, \sigma) \leq \lambda \, {\mathcal E}_1\star {\mathcal T}\, (\mu, \sigma)$.
 
\item  For all $g\in C(X_2)$ and $i\in I$, we have $\beta_1  \big( \int_{X_1} E_1^-\circ T^-g)\, d\mu\big)\leq \frac{1}{\lambda}\beta_2(\int_{X_2} E_2^-(\lambda g) d\nu)$.  
\end{enumerate}
 \end{cor}
\noindent{\bf Proof:} Note that here, we need the formula for $({\mathcal E}\star {\mathcal T})\star (- {\tilde {\cal F}})(\mu, \nu)$. Since $ {\tilde {\cal F}}$ is now a forward convex transfer with Kantorovich operators equal to $\tilde F_i^+(g)=-F_i^-(-g)$, we can apply Part 2) of Proposition \ref{B-F} to ${\cal F}_2=\frac{1}{\lambda} {\tilde{ \mathcal F}}$ and ${\cal F}_1={\mathcal E}\star {\mathcal T}$, which is a $\beta$-backward transfer with Kantorovich operator  $E^-\circ T^-$, to obtain
\[
({\mathcal E}\star {\mathcal T})\star (- {\tilde {\cal F}})(\mu, \nu)=\inf\limits_{g\in C(X_2)}\inf\limits_{j\in J}\left\{-\beta (\int_{X_1} E^-\circ T^-g\, d\mu)+\frac{1}{\lambda}\int_{X_3} F _{j}^-(\lambda g)\, d\nu\right\}.
\]
 A similar argument applies for 2).\\
We now apply the above to the case where $\cal E$ is the backward logarithmic transfer to obtain,
  
\begin{cor} Let  ${\mathcal F}$ be a backward convex transfer on $X_2\times Y_2$ with Kantorovich operators $(F^-_i)_{i\in I}$. 
 Let ${\cal T}$ be a backward linear transfer  on $Y_1\times Y_2$ with Kantorovich operator $T^-$ and $\lambda>0$. Then, for any fixed pair of probability measures $\mu\in {\cal P}( X_1)$ and $\nu\in {\cal P}(X_2 )$, the following are equivalent:
\begin{enumerate}
\item For all $\sigma \in {\cal P}(Y_2)$, we have ${\mathcal F}(\nu, \sigma) \leq \lambda \, {\mathcal H}\star {\mathcal T}\, (\mu, \sigma)$ 
 
 \item  For all $g\in C(X_2)$, we have $\log  \big( \int_{X_1} e^{T^-g}\, d\mu\big)\leq \inf\limits_{i\in I}\frac{1}{\lambda}\int_{X_2} F_i^-(\lambda g) d\nu$.
 \end{enumerate}
\end{cor}

\begin{rmk} \rm An immediate application of (4) in Proposition \ref{forward-back} is the following result in \cite{C-K} 
\begin{equation}
\inf\{{\overline{\mathcal W}}_2(\mu, \sigma)+{\mathcal H}(dx, \sigma); \sigma \in {\cal P}(\R^d)\}=\inf\{-\log \int e^{-f^*}\, dx +\int f\, d\mu; f\in {\cal C}(\R^d)\},
\end{equation} 
where ${\cal C}onv(\R^d)$ is the cone of convex functions on $\R^d$, and ${\overline{\mathcal W}}_2(\mu, \sigma)=-{\mathcal W}_2(\sigma, {\bar \mu})$, the latter being the Brenier transfer of  Example 7.3 and $\bar \mu$ is defined as $\int f(x) d \bar \mu(x)=\int f(-x) d \mu(x)$.  Note that in this case, $T^+f (x)=-f^*(-x)$, $E^-f= e^f$ and $\beta (t)=\log t$, and since $g^{**} \leq g$, 
\begin{eqnarray*}
\inf\{{\overline{\mathcal W}}_2(\mu, \sigma)+{\mathcal H}(dx, \sigma); \sigma \in {\cal P}(\R^d)\}
&=&{\mathcal H}\star (-{\mathcal W}_2)(dx, \bar \mu)\\
&=&\inf\{-\log \int e^{-g}\, dx +\int g^*(x)\, d\mu; g\in C(\R^d)\}\\
&=& \inf\{-\log \int e^{-f^*}\, dx +\int f\, d\mu; f\in {\cal C}onv(\R^d)\}.
\end{eqnarray*} 
What is remarkable in the result of Cordero-Erausquin and Klartag \cite{C-K} is the characterization of those measures $\mu$ (the moment measures) for which there is attainment in both minimization problems. 
\end{rmk}  

\subsection{Maurey-type inequalities}
We are now interested in inequalities of the following type:
For all $\sigma_1\in {\mathcal P}(X_1), \sigma_2\in {\mathcal P}(X_2)$, we have 
\begin{equation}
{\mathcal F}(\sigma_1, \sigma_2) \leq \lambda_1 {\mathcal T}_1\star {\mathcal H}_1(\sigma_1, \mu)+\lambda_2 {\mathcal T}_2\star {\mathcal H}_2( \sigma_2, \nu).
\end{equation}
This will requires a duality formula for the expression
$
\tilde {\mathcal E}_1\star (-{\mathcal T})\star  {\mathcal E}_2,
$
where ${\cal F}$ is a backward convex transfer and $ {\mathcal E}_1$,  ${\mathcal E}_2$ are forward entropic transfers.

\begin{thm} Assume ${\mathcal F}$ is a backward convex transfer on $Y_1\times Y_2$ with Kantorovich operators $(F^-_i)_{i\in I}$, ${\mathcal E}_1$ (resp., ${\mathcal E}_2$) is a forward $\alpha_1$-transfer on $Y_1\times X_1$ (resp., a  forward $\alpha_2$-transfer on $Y_2\times X_2$) with Kantorovich operator  $E_1^+$ (resp., $E_2^+$), then  for any $(\mu, \nu)\in {\cal P}(X_1)\times {\cal P}(X_2)$, we have
\begin{equation}
\tilde {\mathcal E}_1\star (-{\mathcal F}) \star  {\mathcal E}_2 \, (\mu, \nu)=\inf\limits_{i\in I} \inf\limits_{f\in C(X_3)} \left\{\alpha_1 \big(\int_{X_1} E_1^+\circ F_{i}^-f\, d\mu)+ \alpha_2 (\int_{X_2}E_2^+(f)\, d\nu)\right\}.
\end{equation}
 \end{thm}

\noindent{\bf Proof:} If ${\mathcal E}_1$  a forward $\alpha_1$-transfer on $Y_1\times X_1$, then ${\tilde{\mathcal E}}_1$ is a backward $-(\alpha_1^\oplus)^\ominus$-transfer on $X_1\times Y_1$ with Kantorovich operator $\tilde {E_1^-}g=-E_1^+(-g)$. Apply Proposition \ref{back-back} with ${\cal F}_1={\tilde{\mathcal E}}_1$, and ${\cal F}_2={\cal F}$ to get 
\begin{eqnarray*}
\tilde {\mathcal E}_1\star (-{\mathcal F})\, (\mu, \nu)&=&\inf\limits_{f\in C(X_3)}\inf\limits_{i\in I}\left\{(\alpha_1^\oplus)^\ominus \big(\int_{X_1} -E_1^+\circ F_{i}^-f\, d\mu)-\int_{X_3} f\, d\nu \right\}\\
&=&\inf\limits_{f\in C(X_3)}\inf\limits_{i\in I}\left\{\alpha_1 \big(\int_{X_1} E_1^+\circ F_{i}^-f\, d\mu)-\int_{X_3} f\, d\nu \right\}.
\end{eqnarray*}
Write now, 
 \begin{eqnarray*}
\tilde {\mathcal E}_1\star (-{\mathcal F}) \star  {\mathcal E}_2 \, (\mu, \nu)&=&\inf\left\{\tilde {\mathcal E}_1\star (-{\mathcal F}) (\mu, \sigma)+{\mathcal E}_2(\sigma, \nu);\, \sigma \in {\mathcal P}(Y_2)\right\}\\
&=&\inf\limits_{\sigma\in {\mathcal P}(Y_2)} \inf\limits_{f\in C(X_3)}\inf\limits_{i\in I}\left\{\alpha_1 \big(\int_{X_1} E_1^+\circ F_{i}^-f\, d\mu)-\int_{X_3} f\, d\sigma+{\mathcal E}_2(\sigma, \nu)\right\}\\
&=& \inf\limits_{i\in I}\inf\limits_{f\in C(X_3)} \left\{\alpha_1 \big(\int_{X_1} E_1^+\circ F_{i}^-f\, d\mu)-\sup\limits_{\sigma \in {\mathcal P}(Y_2)}\{\int_{Y_2} f\, d\sigma-{\mathcal E}_2(\sigma, \nu)\} \right\}\\
&=&\inf\limits_{i\in I} \inf\limits_{f\in C(X_3)} \left\{\alpha_1 \big(\int_{X_1} E_1^+\circ  F_{i}^-f\, d\mu)+ \alpha_2 (\int_{X_2}E_2^+(-f)\, d\nu)\right\}.
\end{eqnarray*}
 \begin{cor} Assume ${\mathcal E}_1$ (resp., ${\mathcal E}_2$) is a forward $\alpha_1$-transfer on $Z_1\times X_1$ (resp., $\alpha_2$-transfer on $Z_2\times X_2$) with Kantorovich operator $E_1^+$ (resp., $E_2^+$). Let ${\mathcal T}_1$ (resp., ${\mathcal T}_2$) be forward linear transfers on $Y_1\times Z_1$ (resp., $Y_2\times Z_2$) with Kantorovich operator $T_1^+$ (resp., $T_2^+$), and let ${\mathcal F}$ be a backward convex transfer on $Y_1 \times Y_2$ with Kantorovich operators $(F_i^-)_i$. Then, for any given $\lambda_1, \lambda_2\in \R^+$ and $(\mu, \nu)\in  {\mathcal P}(X_1) \times {\mathcal P}(X_2)$, the following are equivalent:
\begin{enumerate}
\item For all $\sigma_1\in {\mathcal P}(Y_1), \sigma_2\in {\mathcal P}(Y_2)$, we have 
\begin{equation}
{\mathcal F}(\sigma_1, \sigma_2) \leq \lambda_1 {\mathcal T}_1\star {\mathcal E}_1(\sigma_1, \mu)+\lambda_2 {\mathcal T}_2\star {\mathcal E}_2( \sigma_2, \nu).
\end{equation}
 \item For all $g\in C(Y_2)$ and all $i\in I$, we have 
\begin{equation}
\lambda_1\alpha_1 \big(\int_{X_1} E_1^+\circ T_1^+\circ (\frac{1}{\lambda_1}F_i^-g)\, d\mu\big)+ \lambda_2\alpha_2 (\int_{X_2} E_2^+\circ T_2^+(\frac{-1}{\lambda_2}g)\, d\nu)\geq 0.
\end{equation}
\end{enumerate}
\end{cor} 
\noindent {\bf Proof:} It suffices to apply the above with the forward $\lambda_i\alpha_i$-transfers ${\cal F}_i:=\lambda_i {\mathcal T}_i\star {\mathcal E}_i$, whose Kantorovich operators are $F_i(g)=E_i^+\circ T_i^+(\frac{g}{\lambda_i})$ for $i=1,2$.\\

By applying the above to ${\cal E}_i(\mu, \nu)=:\cal{H}$ the forward logarithmic entropy where $\alpha_i(t)=-\log(-t)$ and Kantorovich operator $E^+f=e^{-f}$, we get the following extension of a celebrated result of Maurey \cite{Mau}.

\begin{cor} Assume ${\mathcal F}$ is a backward convex transfer on $Y_1\times Y_2$ with Kantorovich operators $(F^-_i)_{i\in I}$, and  
let ${\mathcal T}_1$ (resp., ${\mathcal T}_2$) be forward linear transfer on $Y_1\times X_1$ (resp., $Y_2\times X_2$) with Kantorovich operator  $T_1^+$ (resp., $T_2^+$), then  for any given $\lambda_1, \lambda_2\in \R^+$ and $(\mu, \nu)\in {\cal P}(X_1)\times {\cal P}(X_2)$, the following are equivalent:
\begin{enumerate}
\item For all $\sigma_1\in {\mathcal P}(X_1), \sigma_2\in {\mathcal P}(X_2)$, we have 
\begin{equation}
{\mathcal F}(\sigma_1, \sigma_2) \leq \lambda_1 {\mathcal T}_1\star {\mathcal H}(\sigma_1, \mu)+\lambda_2 {\mathcal T}_2\star {\mathcal H}( \sigma_2, \nu).
\end{equation}

\item For all $g\in C(Y_2)$ and all $i\in I$, we have 
\begin{equation}
(\int_{X_1} e^{-T_1^+\circ \frac{1}{\lambda_1}F_i^-g}\, d\mu\big)^{\lambda_1} (\int_{X_2}e^{-T_2^+(\frac{1}{-\lambda_2}g)}\, d\nu)^{\lambda_2} \leq 1.
\end{equation}
\end{enumerate}
If ${\mathcal T}_1={\mathcal T}_2$ are the identity transfer, then the above is equivalent to saying that for all $g\in C(Y_2)$ and all $i\in I$, we have 
\begin{equation}
(\int_{X_1} e^{\frac{-1}{\lambda_1}F_i^-g}\, d\mu\big)^{\lambda_1} (\int_{X_2}e^{\frac{1}{\lambda_2}g}\, d\nu)^{\lambda_2} \leq 1.
\end{equation}

\end{cor}

\section{Effective transfers and a general discrete weak KAM theory}
 
 Let $X$ be a compact metric space, and let ${\cal T}$ be a backward linear transfer on $\mcal{P}(X) \times \mcal{P}(X)$ with $T$ as a corresponding backward Kantorovich operator.   We will be looking for fixed points of $T$, which correspond to Fathi's notion of weak KAM solutions (\cite{Fa}, \cite{B-B1}) in the case of a transfer induced by a mass transport corresponding to a cost induced by the generating function of a Lagrangian (Example 7.4). 
 
 \begin{thm}\label{WKS} Suppose $\T$ is a backward linear transfer on $\mcal{P}(X) \times \mcal{P}(X)$ that is weak$^*$-continuous on ${\cal M}(X)$, and let $T$ be the corresponding backward Kantorovich operator that maps $C(X)$ into $C(X)$. Then, there exists a constant $ c = c(\T)\in \R$, a backward linear transfer $\T_\infty$ on $\mcal{P}(X) \times \mcal{P}(X)$, and a corresponding Kantorovich operator $T^\infty$ such that 
 \begin{enumerate}
 \item  For all $\mu$ and $\nu$ in ${\mathcal P}(X)$, we have $\T_\infty(\mu,\nu) = (\T + c) \star \T_\infty(\mu,\nu)$ \, and \, $\T_\infty(\mu,\nu) = \T_\infty \star \T_\infty(\mu,\nu)$. 
 \item  For all $f \in C(X)$, we have 
$T^\infty \circ T^\infty f = T^\infty f$. 
\item  The constant $-c(\T)=\inf\{{\cal T}(\mu, \mu); \mu\in {\cal P}(X)\}$
 is attained by a probability measure $\bar{\mu}$ in the set $ {\cal C} := \{\mu\in \mcal{P}(X)\,;\, \T_\infty(\mu, \mu) = 0\}$ such that 
\begin{equation}
(\bar{\mu}, \bar{\mu}) \in \mathcal{D}:= \{ (\mu,\nu) \in \mathcal{P}(X)\times\mathcal{P}(X)\,:\, \T(\mu,\nu) + \T_\infty(\nu,\mu) = -c\}.
\end{equation}
\item For every $f\in C(X)$, the function $u=T^\infty f$ is a solution for the equation 
\begin{equation}\label{ws}
\hbox{$T^nu-nc=u$ for every $n \in \N$,}
\end{equation}
and moreover, every solution of (\ref{ws}) is of this form.
\end{enumerate}
  \end{thm}
  By analogy with the weak KAM theory of Fathi and Mather, we shall say that $\T_\infty$ (resp., $T^\infty$) is {\em the effective transfer} or the {\em generalized Peierls barrier} (resp., {\em effective Kantorovich operator}) associated to $\T$. The constant $c(\T)$ is {\rm the Man\'e critical value}, while the {\em Mather measures} are those $\bar \mu$ where the infimum in (3) is attained. The functions $u$ solving equation  (\ref{ws}) will be called {\em weak KAM solutions} for $\T$. 
   The set ${\mathcal A}$ is the analogue of the {\em projected Aubry set}, and $\mathcal{D}$ can be seen as a {\em generalized Aubry set} \cite{Fa}. 
   
 
 \begin{lem} \label{T_n_Bound}For each $n\in \N$, Let ${\cal T}_n= {\cal T}\star {\cal T}\star ....\star {\cal T}$ be the transfer obtained from $\T$ by iterating its convolution $n$-times. Then,
 \begin{enumerate}
 
 \item For all $\mu, \nu \in {\cal P}(X)$, we have
$ {\cal T}_n (\mu, \nu)=\sup\big\{\int_{X}g(y)\, d\nu(y)-\int_{X}T^ng(x);\,  g\in C(X) \big\}.$

 \item The sequence of transfers $(\T_n)_n$ is equicontinuous, and there exists a positive constant $C>0$ and a number $c \in \R$ such that 
  \begin{equation}
\hbox{$| {\cal T}_n (\mu, \nu) -c n| \leq C$ \quad for all $\mu, \nu \in {\cal P}(X)$ and $n\in \N$.}
 \end{equation}
 
 \end{enumerate}
 
 \end{lem}
 \noindent{\bf Proof:} 1) is  immediate from Proposition \ref{inf.tens}.\\
For 2), we follow an argument of  Bernard-Buffoni \cite{B-B1}. Since $\T$ is weak$^*$ continuous and $X$ is compact, there exists a modulus of continuity $\delta: [0,\infty) \to [0,\infty)$, with $\lim_{t \to 0}\delta(t) = \delta(0) = 0$ such that 
\begin{equation}
|\T(\mu,\nu) - \T(\mu',\nu')| \leq \delta(W_2(\mu,\mu') + W_2(\nu,\nu')),\quad \text{for all } \mu, \mu', \nu, \nu' \in \mathcal{P}(X).
\end{equation}
Since for all $\sigma_1,\ldots,\sigma_{n-1} \in \mathcal{P}(X)$, the map
\begin{equation*}
(\mu,\nu) \mapsto \T(\mu,\sigma_1) + \T(\sigma_1,\sigma_2) + \ldots + \T(\sigma_{n-1},\nu)
\end{equation*}
also has the same modulus of continuity $\delta$, then $\T_n$ does too, as the infimum of functions with the same modulus of continuity, also has the same modulus of continuity.
 Define the sequences 
 \begin{align*}
 M_n &:= \max\left\{\T_n(\mu,\nu)\,;\, \mu,\nu \in \mathcal{P}(X)\right\}\\
 m_n &:= \min\left\{\T_n(\mu,\nu)\,;\, \mu,\nu \in \mathcal{P}(X)\right\}.
 \end{align*}
 Then $M_n$ is sub-additive, while $m_n$ is super-additive. Indeed, there exists $\mu,\nu \in \mathcal{P}(X)$ such that, for all $\sigma \in \mathcal{P}(X)$,
 \begin{equation*}
 M_{n+k} = \T_{n+k}(\mu,\nu) \leq \T_{n}(\mu,\sigma) + \T_{k}(\sigma,\nu) \leq M_n + M_k.
 \end{equation*}
The argument for $m_n$ is the same with reverse inequalities.\\
By the equicontinuity of $\T_n$, the difference $M_n - m_n$ is bounded above by a constant $C$ independent of $n$. Therefore it is well known that $\frac{M_n}{n}$ decreases to $M := \inf_{n}\frac{M_n}{n}$, while $\frac{m_n}{n}$ increases to $m := \sup_{n}\frac{m_n}{n}$. Then the string of inequalities
\begin{equation*}
nM - C \leq M_n - C \leq m_n \leq \T_n(\mu,\nu) \leq M_n \leq m_n + C \leq nm + C
\end{equation*}
implies $M \leq m + \frac{2C}{n}$ for all $n$; hence $c:=m = M$.

\begin{prop}\label{0case} Under the condition that $c=0$, there exists a backward linear transfer ${\cal T}_\infty$ on $\mcal{P}(X) \times \mcal{P}(X)$ with associated Kantorovich operator $T^\infty$ satisfying the following properties.

\begin{enumerate}

\item For every  $f \in C(X)$, we have
\hbox{$T \circ T^\infty f = T^\infty f$ \quad and \quad $T^\infty \circ T^\infty f = T^\infty f$.} The operator $T^\infty$ is given explicitly by 
 \begin{equation}\label{explicit.form.Tinfty}
 T^\infty f(x) = \lim_{n \to \infty}T^n  \circ \limsup_{m \to \infty}T^m f(x),
\end{equation}
 and $\T_\infty(\mu,\nu) = \sup\lf\{\int_{X}f\d\nu - \int_{X}T^\infty f\d\mu\rt\}$.
Moreover, the fixed points of $T$ and the fixed points of $T_\infty$ are the same. 

\item For all $\mu, \nu \in {\cal P}(X)$ and $n \in \N$, we have 
\begin{equation}
\hbox{$\T_\infty(\mu,\nu) = \T_n * \T_\infty(\mu,\nu)$ \quad and \quad $\T_\infty(\mu,\nu) = \T_\infty * \T_\infty(\mu,\nu)$.}
\end{equation}

\item $\sup\lf\{\int_{X}T^\infty f\d (\nu-\mu)\,;\, f \in C(X)\rt\} \leq \T_\infty(\mu,\nu) \leq \liminf_{n \to \infty}\T_n(\mu,\nu)$.

\item The set ${\cal A}:=\{\mu \in {\cal P}(X); {\cal T}_\infty(\mu, \mu)=0\}$ is non-empty and for every 
$\mu, \nu \in {\cal P}(X)$, we have 
\begin{equation}
{\cal T}_\infty(\mu, \nu)=\inf\{ {\cal T}_\infty(\mu, \sigma)+{\cal T}_\infty(\sigma, \nu), \sigma \in {\cal A}\},
\end{equation}
and the infimum on ${\cal A}$ is attained. 
\item We have 
\begin{equation}
\inf\{{\cal T}(\mu, \mu); \mu\in {\cal P}(X)\}=0,
\end{equation}
and the infimum is attained by a measure $\bar{\mu} \in \mathcal{A}$ such that 
\begin{equation}
(\bar{\mu}, \bar{\mu}) \in \mathcal{D}:= \{ (\mu,\nu) \in \mathcal{P}(X)\times\mathcal{P}(X)\,:\, \T(\mu,\nu) + \T_\infty(\nu,\mu) = 0\}.
\end{equation}
 \end{enumerate}
\end{prop} 

\noindent{\bf Proof:} 
(1)  By Lemma \ref{T_n_Bound}, $\T_n$ is equicontinuous with modulus of continuity $\delta$. Therefore, for each $f \in C(X)$, $x \mapsto T^n f(x)$ is uniformly continuous with the same modulus of continuity $\delta$. Moreover, $f \mapsto T^n f$ is $1$-Lipschitz with respect to the $\sup$-norm. Indeed, 
\begin{align*}
T^n f(x) &=  \sup_{\nu \in \mathcal{P}(X)}\left\{\int_{X}f\,d\nu - \T_n(\delta_x,\nu)\right\}\\
&=  \sup_{\nu \in \mathcal{P}(X)}\left\{\int_{X}f\,d\nu - \T_n(\delta_y,\nu) + \T_n(\delta_y,\nu) - \T_n(\delta_x,\nu)\right\}\\
&\leq \sup_{\nu \in \mathcal{P}(X)}\left\{\int_{X}f\,d\nu - \T_n(\delta_y,\nu)+ \delta(d(x,y))\right\}\\
&= T^nf (y) + \delta(d(x,y)).
\end{align*}
Interchanging $x$ and $y$ we conclude the continuity. For $1$-Lipschitz,
\begin{align*}
T^n f(x) &=  \sup_{\nu \in \mathcal{P}(X)}\left\{\int_{X}f\d\nu - \T_n(\delta_x,\nu)\right\}\\
&\leq  \sup_{\nu \in \mathcal{P}(X)}\left\{\int_{X}g\d\nu - \T_n(\delta_x,\nu)\right\} + \|f-g\|_{\infty}\\
&= T^ng (x) + \|f-g\|_{\infty}.
\end{align*}
Interchanging $f$ and $g$ we conclude.\\
Define $T_1^\infty f(x) := \limsup_{n \to \infty} T^n f(x)$. The assumption $c = 0$ ensures that $T_1^\infty f(x) < \infty$, and $T_1^\infty$ satisfies the same properties above as $T$. We have the following monotonicity:
\begin{equation}
T\circ T_1^\infty f(x) \geq T_1^\infty f(x),\quad \text{for all $f \in C(X)$ and all $x \in X$.}
\end{equation}
Indeed, since $\left\{\sup_{k \geq n} T^kf\right\}_{n}$ is a sequence of continuous functions which pointwise decreases monotonically to the continuous function $T_1^\infty f$, the sequence $\sup_{k \geq n} T^kf$ converges uniformly to $T_1^\infty f$. We conclude by the Lipschitz property of $T$.\\
Therefore, for each $f \in C(X)$ and $x \in X$, $\{T^n\circ T_1^\infty f(x)\}_{n}$ is a monotone sequence. The assumption $c = 0$ implies it is uniformly bounded in $n$, hence the pointwise limit exists and is finite, and we may define an operator $T^\infty$ by
\begin{equation}
T^\infty f(x) := \lim_{n \to \infty} T^n \circ T_1^\infty f(x).
\end{equation}
Note that again $T^\infty$ satisfies the same properties as $T$; in particular it is a convex operator. Moreover, $x \mapsto T^\infty f(x)$ is continuous, and therefore, the convergence is uniform. Then, by the Lipschitz property of $T$, we get that $T\circ (T^n \circ T_1^\infty) f$ converges uniformly to $T\circ T^\infty f$. In other words,
\begin{equation}
T^\infty f = \lim_{n \to \infty} T^{n}\circ T_1^\infty f = \lim_{n \to \infty} T^{n+1} \circ T_1^\infty f = T\circ T^\infty f.
\end{equation}
Finally, suppose $f \in C(X)$ is a fixed point of $T$. Then $T^n f = f$, so $T_1^\infty f = f$, and consequently $T^n \circ T_1^\infty f = f$. Letting $n \to \infty$, we get $T^\infty f = f$.\\
Conversely, suppose $f$ is a fixed point of $T^\infty$. Since $T^\infty f$ is a fixed point of $T$ from above, we get that $T f = f$.

(2) Define
\begin{equation}
\T_\infty(\mu,\nu) := \sup\left\{\int_{X}f\,d\nu - \int_{X}T^\infty f\,d\mu\,;\, f \in C(X)\right\}.
\end{equation}
Since $T^\infty$ is a convex operator, $\T_\infty$ is a backward linear transfer. From 1), we get immediately the conclusion of 2), by Proposition \ref{inf.tens}.

(3) Note that $T^\infty f \geq T_1^\infty f = \limsup_{n} T^n f$, so
\begin{align*}
\T_\infty(\mu,\nu) &= \sup_{f \in C(X)}\left\{\int_{X}f\d\nu - \int_{X}T^\infty f\d\mu\right\}\\
&\leq \sup_{f \in C(X)}\left\{\int_{X}f\d\nu - \int_{X}\limsup_{n}T^n f\d\mu\right\}\\
&\leq \sup_{f \in C(X)}\liminf_{n}\left\{\int_{X}f\d\nu - \int_{X}T^n f\d\mu\right\}\quad \text{(Fatou)}\\
&\leq \liminf_{n}\sup_{f \in C(X)}\left\{\int_{X}f\d\nu - \int_{X}T^n f\d\mu\right\}\\
&= \liminf_{n}\T_n(\mu,\nu).
\end{align*}
On the other hand, from $T^\infty \circ T^\infty f = T^\infty f$,
\begin{align*}
\T_\infty(\mu,\nu) &= \sup\lf\{\int_{X}f\d\nu - \int_{X}T^\infty f\d\mu\,;\, f \in C(X)\rt\}\\
&\geq \sup\lf\{\int_{X}T^\infty f\d(\nu-\mu)\,;\, f \in C(X)\rt\}.
\end{align*}
(4) This argument is a minor modification of the one given in \cite{B-B1}, to account for the fact that $\T_\infty$ is in general only weak$^*$ lower semi-continuous and not weak$^*$ continuous. Fix $\mu, \nu \in \mathcal{P}(X)$. By 2), there exists $\sigma_1 \in \mathcal{P}(X)$ such that 
\begin{equation*}
\T_\infty(\mu,\nu) = \T_\infty(\mu,\sigma_1) + \T_\infty(\sigma_1, \nu).
\end{equation*}
Similarly, there exists a $\sigma_2$ such that 
\begin{equation*}
\T_\infty(\sigma_1, \nu) = \T_\infty(\sigma_1, \sigma_2)+ \T_\infty(\sigma_2, \nu).
\end{equation*}
Combining the above two equalities, we obtain
\begin{equation*}
\T_\infty(\mu,\nu) = \T_\infty(\mu,\sigma_1) + \T_\infty(\sigma_1, \sigma_2)+ \T_\infty(\sigma_2, \nu).
\end{equation*}
Note also that
\begin{equation}\label{split}
\T_\infty(\mu,\sigma_1) + \T_\infty(\sigma_1, \sigma_2) = \T_\infty(\mu,\sigma_2).
\end{equation}
This follows from 
\begin{align*}
\T_\infty(\mu,\nu) &= \T_\infty(\mu,\sigma_1) + \T_\infty(\sigma_1, \sigma_2)+ \T_\infty(\sigma_2, \nu)\\
&\geq \T_\infty*\T_\infty(\mu, \sigma_2) + \T_\infty(\sigma_2, \nu)\\
&= \T_\infty(\mu, \sigma_2) + \T_\infty(\sigma_2, \nu)\\
&\geq \T_\infty*\T(\mu,\nu)\\
&= \T_\infty(\mu,\nu).
\end{align*}
Hence all the inequalities are equalities; in particular (\ref{split}).\\
After $k$ times we have
\begin{equation*}
\T_\infty(\mu,\nu) = \sum_{i = 0}^{k}\T_\infty(\sigma_{i},\sigma_{i+1}) 
\end{equation*}
where $\sigma_0 := \mu$ and $\sigma_{k+1} := \nu$. This inductively generates a sequence $\{\sigma_{k}\}$ with the property 
\begin{equation*}
\sum_{i = \ell}^{m} \T_\infty(\sigma_{i}, \sigma_{i+1}) =  \T_\infty(\sigma_{\ell}, \sigma_{m+1})
\end{equation*}
whenever $0 \leq \ell < m \leq k$. In particular, for any subsequence $\sigma_{k_j}$, we have
\begin{equation}\label{sebseq}
\T(\mu,\sigma_{k_1}) + \sum_{j = 1}^{m} \T_\infty(\sigma_{k_j}, \sigma_{k_{j+1}}) + \T_\infty(\sigma_{k_{m+1}}, \nu) =  \T_\infty(\mu, \nu).
\end{equation}
Extract a weak$^*$ convergent subsequence $\{\sigma_{k_j}\}$ to some $\bar{\sigma} \in \mathcal{P}(X)$. By weak$^*$ l.s.c. of $\T_\infty$, we have
\begin{equation*}
\liminf_{j}\T_\infty(\sigma_{k_j}, \sigma_{k_{j+1}}) \geq \T_\infty(\bar{\sigma},\bar{\sigma}).
\end{equation*}
In particular, given $\epsilon > 0$, for all but finitely many $j$,
\begin{equation}\label{biggerthan}
\T_\infty(\sigma_{k_j}, \sigma_{k_{j+1}}) \geq \T_\infty(\bar{\sigma},\bar{\sigma}) - \epsilon.
\end{equation}
Therefore, by refining to a further (non-relabeled) subsequence if necessary, we obtain a subsequence $\{\sigma_{k_j}\}$ satisfying (\ref{biggerthan}) for all $j$. By further refinement, we may also assume,
\begin{equation}\label{boundarybigger}
\T_\infty(\mu, \sigma_{k_1}) \geq \T_\infty(\mu, \bar{\sigma}) - \epsilon.
\end{equation}
Therefore, by refining to a further (non-relabeled) subsequence if necessary, we obtain a subsequence $\{\sigma_{k_j}\}$ with properties (\ref{sebseq}), (\ref{biggerthan}), and (\ref{boundarybigger}).

Moreover, for all $m$ large enough (depending on $\epsilon$), we have
\begin{equation}\label{tailend}
\T_\infty(\sigma_{k_{m+1}}, \nu) \geq \T_\infty(\bar{\sigma}, \nu) - \epsilon
\end{equation}
Applying the inequalities of (\ref{biggerthan}), (\ref{boundarybigger}), and (\ref{tailend}), to (\ref{sebseq}), we obtain
\begin{equation*}
\T_\infty(\mu,\nu) \geq \T_\infty(\mu,\bar{\sigma}) + m\T_\infty(\bar{\sigma},\bar{\sigma}) + \T_\infty(\bar{\sigma},\nu)- (m+2)\epsilon
\end{equation*}
for all large enough $m$. From the fact that $\T_\infty = \T_\infty * \T_\infty$, the above inequality is only possible if 
\begin{equation*}
\T_\infty(\bar{\sigma},\bar{\sigma}) \leq \frac{m+2}{m}\epsilon \leq 2\epsilon.
\end{equation*}
As $\epsilon$ is arbitrary, we obtain $\T_\infty(\bar{\sigma},\bar{\sigma}) \leq 0$, and consequently $\T_\infty(\bar{\sigma},\bar{\sigma}) = 0$ (the reverse inequality following from $\T_\infty = \T_\infty * \T_\infty$). 

Finally, we note that $\T_\infty (\mu,\nu) = \T_\infty(\mu,\sigma_{k_j}) + \T_\infty(\sigma_{k_j},\nu)$ for all $j$, so at the $\liminf$, we find
\begin{equation*}
\T_\infty(\mu,\nu) \geq \T_\infty(\mu,\bar{\sigma}) + \T_\infty(\bar{\sigma},\nu).
\end{equation*}
The reverse inequality is immediate from $\T_\infty = \T_\infty * \T_\infty$.

(5)  
By the assumption $c = 0$, we have $\T(\mu,\mu) \geq 0$ for all $\mu$. Therefore it suffices to find $\bar{\mu}$ such that $\T(\bar{\mu}, \bar{\mu}) = 0$. We may construct a sequence $(\mu_k) \subset \mathcal{A}$ such that $(\mu_k,\mu_{k+1}) \in \mathcal{D}$. The set $\mathcal{D}$ is convex by convexity of both $\T$ and $\T_\infty$. Therefore, the Cesaro averages $(\frac{1}{n}\sum_{k = 1}^{n}\mu_k,\frac{1}{n}\sum_{k = 1}^{n}\mu_{k+1}) $ belong to $\mathcal{D}$.\\ 

Denoting $\nu_n := \frac{1}{n}\sum_{k = 1}^{n}\mu_k$, we have
$(\nu_n, \nu_n + \frac{1}{n}(\mu_{n+1} - \mu_1)) \in \mathcal{D}.$
 Extracting a weak$^*$ convergent subsequence $\nu_{n_j}$ converging to some $\bar{\mu} \in \mathcal{A}$ (since, in particular $\nu_n \in \mathcal{A}$ and $\mathcal{A}$ is weak$^*$ closed), we use the fact that $(\nu_n, \nu_n + \frac{1}{n}(\mu_{n+1} - \mu_1)) \in \mathcal{D}$ and weak$^*$ continuity of $\T$ (resp. weak$^*$ l.s.c. of $\T_\infty$) to find
\begin{equation*}
\T(\bar{\mu},\bar{\mu}) \leq -\T_\infty (\bar{\mu}, \bar{\mu}) = 0,
\end{equation*}
which concludes the proof.

Theorem \ref{WKS} now follows from Proposition \ref{0case} by replacing the transfer $\T$ by $\T +c$ and $T$ by $T -c$.

\begin{rmk} {\bf The case of a cost minimizing mass transport:} \rm
Consider the setting of a regular mass transport, with a continuous cost $A$ on $X \times X$.
\eqs{
\T(\mu,\nu) := \inf_{\pi\in \mcal{K}(\mu,\nu)}\int_{X \times X} A(x,y)\d\pi(x,y).
}
Then, under the assumption that $c = 0$, a computation yields that
\eqs{
T^\infty f(x) = \sup_{y}\{f(y) - A_\infty(x,y)\}\quad\text{where}\quad A_\infty(x,y) := \liminf_{n \to \infty}A_n(x,y),
}
and the generalized Peierl's barrier is an optimal mass transport with cost $A_\infty(x,y)$,
\eqs{
\T_\infty(\mu,\nu) = \inf_{\pi \in \Gamma(\mu,\nu)}\int A_\infty(x,y)\d\pi(x,y).
}
Therefore, the objects $\T_\infty$, $T^\infty$, above, reduce to those studied by Bernard-Buffoni \cite{B-B2}.  

\end{rmk}

\begin{rmk} {\bf The case of Pushforward Transfers}\rm

Recall the pushforward transfer of Example 1: For $\sigma: X \to X$ a continuous map,
 \begin{equation}
{\mathcal I}(\mu, \nu)=\left\{ \begin{array}{llll}
0 \quad &\hbox{if $\sigma_\#\mu=\nu $}\\
+\infty \quad &\hbox{\rm otherwise.}
\end{array} \right.
\end{equation}
This transfer is not weak$^*$ continuous so strictly speaking the above theorem does not apply, however the operator $T^\infty$ as defined in (\ref{explicit.form.Tinfty}) can still make sense, and it is of the form $T^\infty f = f\circ \sigma^\infty$, where $\sigma^\infty(x) := \limsup_{m \to \infty}\sigma^m(x)$. 

If one takes $X = [0,1] \subset \R$ and $\sigma(x) = x^2$, then $\sigma^\infty(x) = 0$, if $x \in [0,1)$ and $1$ if $x = 1$. If $\sigma(x) = 1-x$, then $\sigma^\infty (x) = \max\{x, 1-x\}$. In both examples, it is easy to check that $T^\infty f$ is a fixed point of $T$. Note that for $\sigma(x) = x^2$, $\sigma^\infty$ is not continuous, yet $T^\infty$ is \textit{not} satisfies the other convexity properties but not on spaces of continuous functions. This points towards a possible extension of the theory beyond this class and to other spaces of measurable functions
The role of the weak$^*$ continuity assumption on $\T$ was to ensure that $T^\infty$ is a backward Kantorovich operator in the definition we have adopted throughout this paper.
\end{rmk}

\end{document}